\documentclass[11pt]{amsart}
\usepackage{amsmath}
\usepackage{amscd}
\usepackage{amssymb}
\usepackage{graphicx}
\usepackage{graphics}
\usepackage{latexsym}
\usepackage[hyphens]{url}

\input xy
\xyoption{all}

\newtheorem{theorem}{Theorem}[section]
\newtheorem{lemma}[theorem]{Lemma}

\theoremstyle{definition}
\newtheorem{definition}[theorem]{Definition}
\newtheorem{remark}[theorem]{Remark}

\begin{document}
\title{Measuring Congressional district meandering}
\author{Eion Blanchard}
\email{eionmb2@illinois.edu}
\author{Kevin Knudson}
\email{kknudson@ufl.edu }
\address{Department of Mathematics, University of Florida, P.O.~Box 118105, Gainesville, FL 32611-8105}

\date{\today}

\begin{abstract} In recent decades, state legislatures have often drawn U.S. Congressional voting districts that look---to the human eye---to be rather twisted. In this paper, we propose a method to measure how much districts ``meander'' via a computation of the medial axis of the region. We then compare this to the medial axis of the convex hull of the district to obtain the {\em medial-hull ratio}: a dimensionless quantity that captures the district's irregularity. We compute this quantity for many example Congressional districts.
\end{abstract}

\maketitle

\section{Introduction}\label{sec:intro}
Since 1812, when Massachusetts Governor Elbridge Gerry signed off on a state redistricting plan featuring a salamander-shaped region clearly designed to favor his party's chances in future elections, politicians have angled to draw lines on the map in self-serving ways. In the intervening two centuries, state legislatures have often pushed this practice to its limits, creating increasingly precise and complicated legislative districts in order to cement their power. 

The last couple of decades have witnessed attempts by mathematicians, statisticians, and political scientists to quantify the extent to which a district has been gerrymandered. Various metrics have been developed: the Polsby-Popper score, which is essentially the isoperimetric ratio of a region's area to the square of its perimeter; the convex hull measure of the ratio of a region's area to that of its convex hull; and the efficiency gap, which seeks to quantify the proportion of wasted votes in each district. The first two of these are geometric in nature and are measures of the ``compactness'' of a district. While this concept is a necessary feature of Congressional districts, as outlined in the Voting Rights Act, no court has given an adequate definition of what it means for a region to be compact. Luckily, mathematicians have a good understanding of this idea.

\begin{figure}
\centerline{\includegraphics[width=4in]{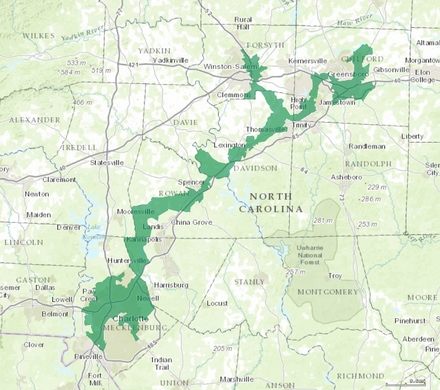}}
\caption{\label{nc12} North Carolina district 12}
\end{figure}

Some districts have drawn significant attention for how much they twist and meander through a state's area. One famous example in just the past five years is the former North Carolina district 12 (Figure \ref{nc12}), which stretched from Durham in the northeast to Charlotte in the southwest and reached up in the northwest to grab a piece of Winston-Salem. To be sure, this district looks odd, and whenever a map is drawn this way, it attracts a crowd of observers. 

There have been few explicit attempts to measure how much a district ``meanders'' in this way. One possible approach was outlined by Dr.~Lisa Handley in a footnote of an expert witness report for a 1994 court case \cite{handley}; we describe this in Section \ref{sec:meander}. Her metric does measure something related to how twisted a map may be, but the court was unconvinced in the case. 

In this paper we propose a metric, the {\em medial-hull ratio}, as a possible solution to this question. The idea is simple: given a polygonal region in the plane we consider its medial axis, which can be thought of as a type of spine for the shape. For a region that meanders significantly we should expect a long medial axis, even if the region is contained in a relatively small area (think of a space-filling curve). We then compare the length of the medial axis to that of the medial axis of the convex hull of the region to obtain our metric. A large value of this ratio indicates that the district is wandering around quite a bit inside the convex area containing the region, suggesting that further investigation into the region's boundary may be warranted.

This paper is organized as follows. In Section \ref{sec:meander}, we describe Handley's metric and elucidate some potential problems with it. Section \ref{sec:medial} gives a definition of the medial axis and indicates how it may be computed. We describe our procedure for computing the medial axis of a Congressional district in Section \ref{sec:methods}; this includes a discussion of various technical issues that must be overcome and potential pitfalls in the calculations. Section \ref{sec:computation} contains the results of computational experiments to determine these ratios for a variety of Congressional districts. Particular attention is given to Pennsylvania, whose Supreme Court recently ordered a new map drawn after ruling the existing districting plan violated the state constitution. Our computations show that the new map is fairer under our measure. One might still argue about whether districts have been gerrymandered, but not on mathematical grounds.

\section{Handley's meanderingness test}\label{sec:meander} We describe here the ``meanderingness measure'' as defined in a footnote in Handley's expert witness report in Johnson v.~Miller \cite{handley}. It is meant to measure how much a district twists and turns, and therefore should reflect a lack of regularity of the region. 

Congressional districts are built using census blocks as constructed by the United States Census Bureau. The algorithm begins with the centroids of each census block in the district, which are readily available from the Census Bureau. Given a district $D$, let $B$ be a census block with centroid $P$. Draw a ray due north from $P$ extending beyond the border of the district and find the first intersection point of this ray with the boundary. It may be that the ray hits the boundary in multiple points, but we consider only the intersection point closest to the centroid. Repeat this procedure by moving around the compass in $5$-degree increments to obtain a collection of 72 points around the centroid of the block on the boundary of the district. These points then define a {\em coverage polygon} $C(B)$ whose boundary is constructed by connecting adjacent vertices with straight lines. The idea is that this polygon approximates the portion of the district that can be reached from the centroid via straight lines.

Now for each census block $B$, compute the area $\alpha(C(B))$ and the ratio $\alpha(C(B))/\alpha(D)$ of this area to the total area of the district. Repeat this for each 25th census block in the district; denote the set of blocks in this sample by ${\mathcal B}$. 

\begin{definition}\label{def:meandermeasure} The {\em meanderingness measure} of $D$ is defined to be $$\mu(D) = \max_{B\in{\mathcal B}} \biggl\{\frac{\alpha(C(B))}{\alpha(D)}\biggr\}.$$
\end{definition}

\begin{remark} Note the counterintuitive fact that a high value of $\mu(D)$ represents a district $D$ that ``meanders'' less; that is, a large value indicates that there is a block from which a large portion of the district is reachable by straight lines. So we would say a district meanders more the {\em smaller} its $\mu$ value is.
\end{remark}

One immediately notes several problems with Definition \ref{def:meandermeasure}. First, why compute the coverage polygon for every 25th census block? Handley's rationale is that it is computationally infeasible to compute these for every block in the district. Indeed, in Johnson v.~Miller the district in question contained more than 18,000 census blocks, so this was a valid consideration (especially with 1994 computing power). Still, no mention is made of how this sampling was conducted. ``Every 25th block'' implies a linearly ordered list of districts, presumably from some order imposed by the Census Bureau. How arbitrary is this ordering? It is conceivable that the blocks are numbered in some way so that the collection ${\mathcal B}$ clusters in a small part of the district. One assumes, however, that the sample ${\mathcal B}$ is chosen to guarantee a reasonable spatial distribution.

\begin{figure}
\centerline{\includegraphics[height=3in]{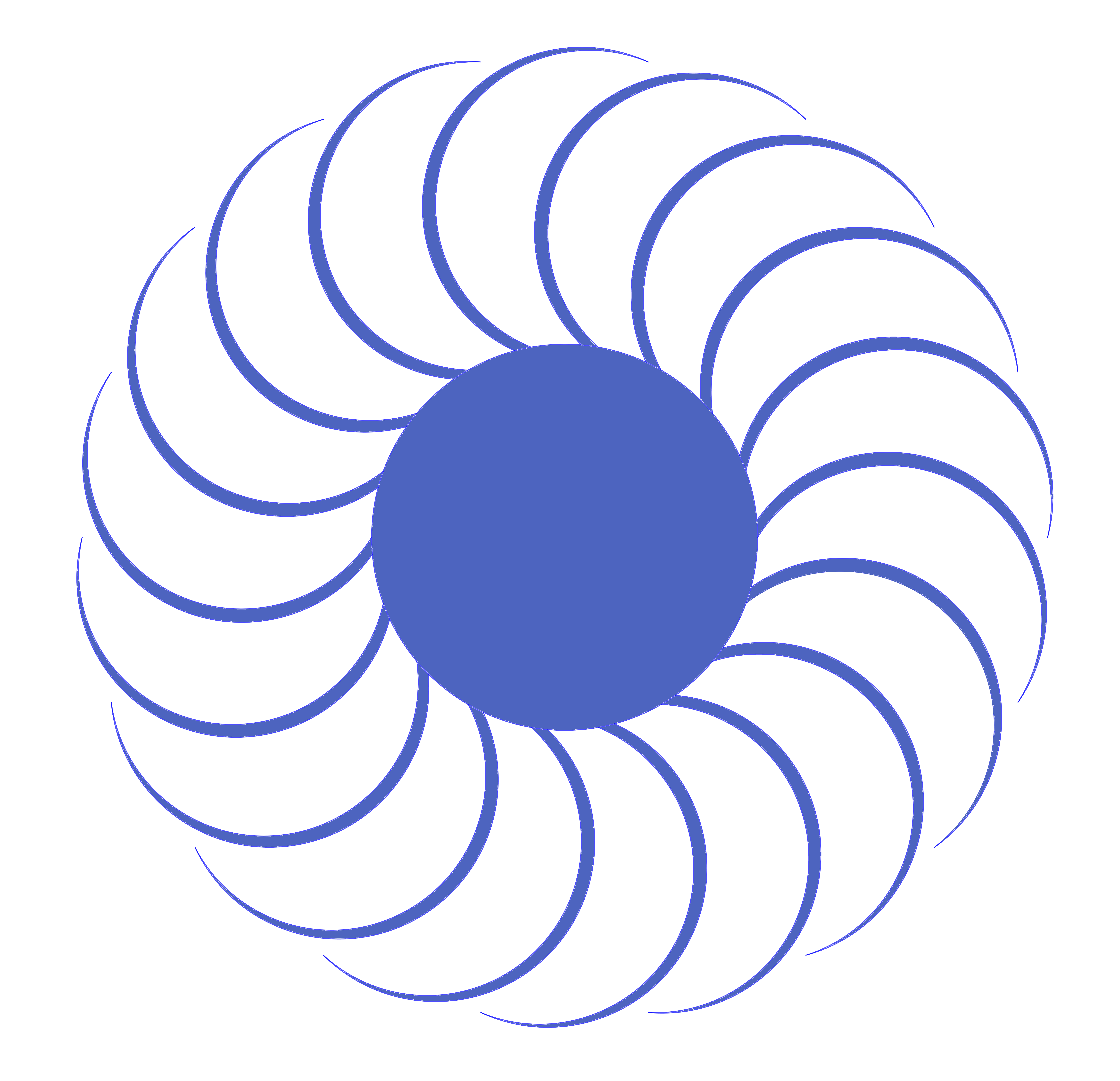}}
\caption{\label{fig:badblock} A district with meanderingness measure $\mu\approx 0.61$.}
\end{figure}

More troubling, however, is the choice of vertices of the coverage polygon of block $B$. Again, for computational purposes, a choice must be made, but consider the district shown in Figure \ref{fig:badblock}. No matter how many census blocks make up the district, a few things are clear. If a block has its centroid in one of the curved arms, then the corresponding coverage polygon will be very small and will therefore not be a candidate to yield the value $\mu(D)$. On the other hand, if a block has its centroid somewhere in the central circular region then the coverage polygon will be a $72$-gon approximating the circle. We may calculate directly that $\mu(D)\approx 0.61$, a relatively large value, yet the district clearly ``meanders.'' Fewer spiral arms leads to a larger value for $\mu(D)$, while still yielding a district that wanders around. While one might object that this is an unrealistic example, it does illustrate the point that such long tentacles will distort this measure. Indeed, the ruling in the case found this metric unconvincing at best, specifically pointing out that long, narrow districts hundreds of miles long would be considered compact under this measure.

\section{The Medial Axis}\label{sec:medial}

\begin{definition}
The {\sl medial axis} of a planar region $R$ is the locus of centers of maximal circles (touching the boundary of $R$ in at least two points) contained in $R$.
\end{definition}

There are a number of equivalent definitions and conceptualizations of the medial axis. For instance, one likens the medial axis to the aftermath of a grassfire: imagine setting a grassy region's boundary ablaze all at once. The quench points---those where the fire, uniformly moving inward, intersects and extinguishes itself---wholly compose the medial axis. Others name the medial axis as a region's skeleton or symmetric axis.

While medial axes may be taken on regions with free boundary, we acknowledge that U.S. Congressional district boundaries are delineated as simple polygons, so our scope in this paper will be limited to the medial axes of these figures.

When we consider the possible types of maximal circles in a simple polygon, we witness those touching two edges, those touching one edge and one vertex, and those touching two vertices. A useful consequence of this is the following lemma, the proof of which may be found in \cite{diffeo}.

\begin{lemma}
For a simple polygon, the medial axis is a union of line segments and portions of parabolic curves. 
\end{lemma}

\subsection{Examples}
\begin{figure}
\centerline{\includegraphics[height=1.5in]{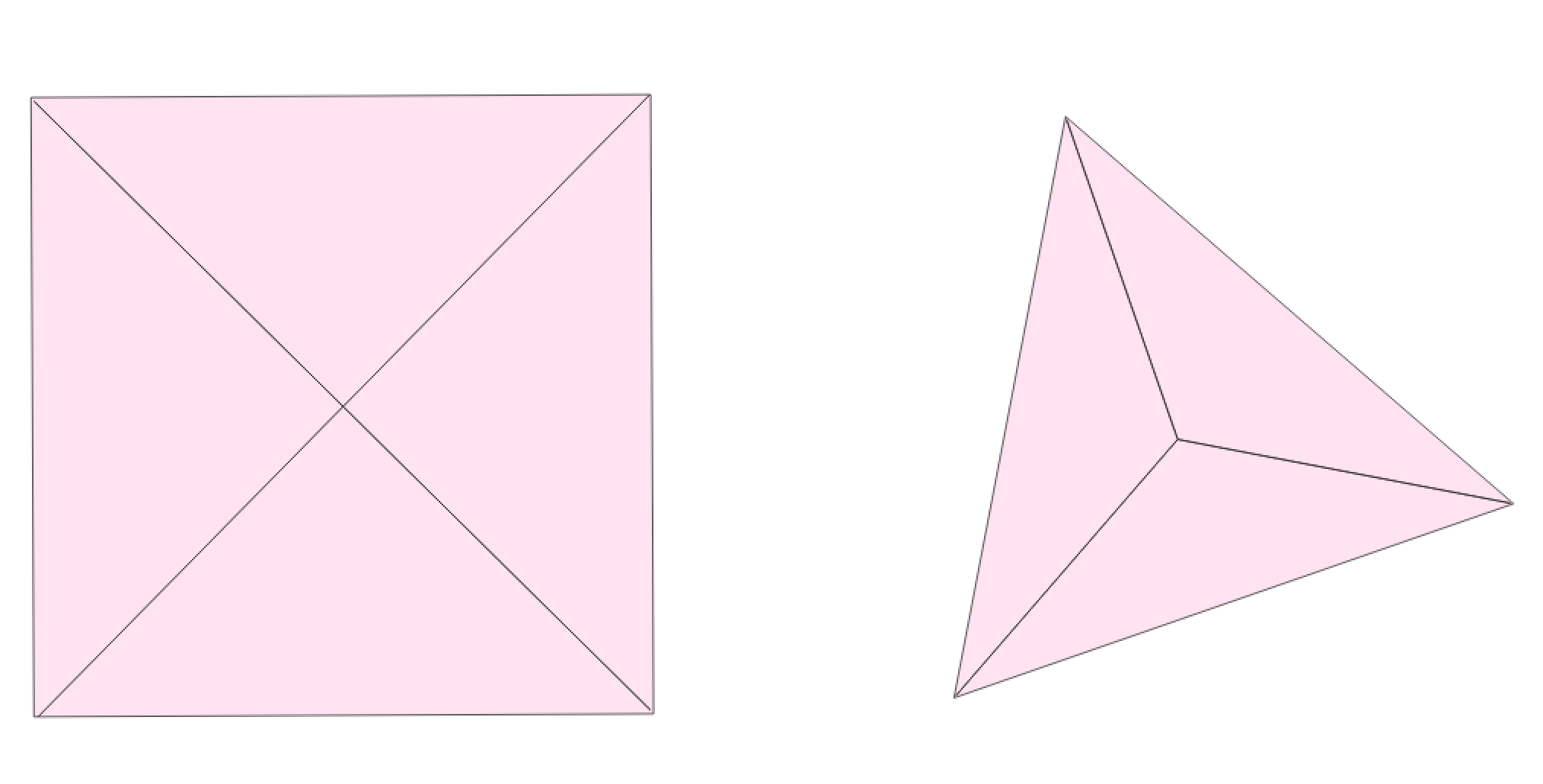}}
\caption{\label{fig:squ_tri} A square and an equilateral triangle with their medial axes: the square's medial axis is simply its diagonal lines, and the triangle's medial axis is the set of line segments touching the centroid and bisecting each interior angle.}
\end{figure}

\begin{figure}
\centerline{\includegraphics[height=1.6in]{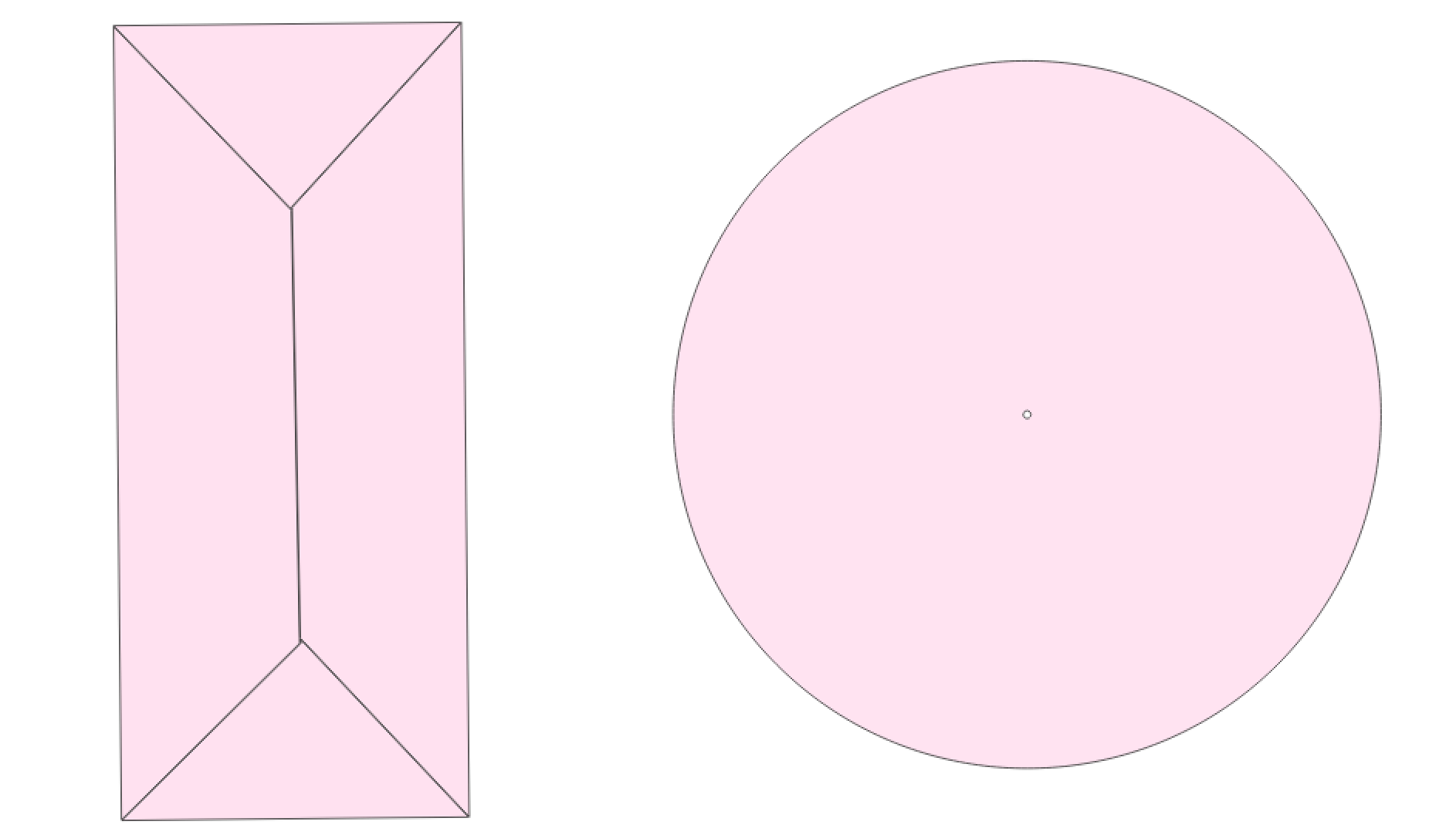}}
\caption{\label{fig:rec_cir} A rectangle and a circle with their medial axes: the rectangle's medial axis is the central line along the longer side, which splits on either end into two line segments so that the interior angles are bisected; and the circle's medial axis is simply its center point.}
\end{figure}

\begin{figure}
\centerline{\includegraphics[height=2.0in]{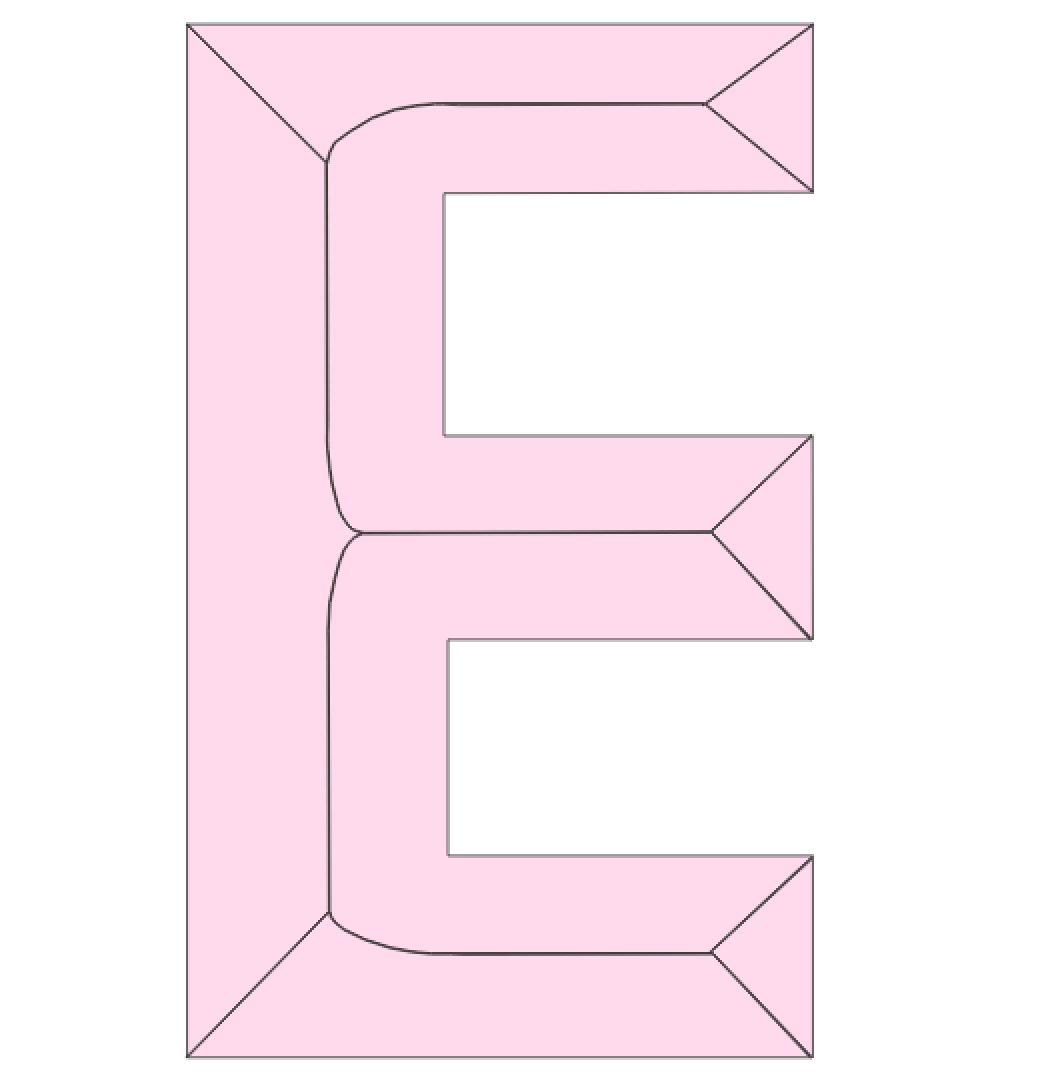}}
\caption{\label{fig:E} A letter `E' with its medial axis: it is composed of rectangle-like axes on the arms, square-like diagonals on the outer corners, and portions of parabolas on the inner corners.}
\end{figure}

\begin{figure}
\centerline{\includegraphics[height=1.8in]{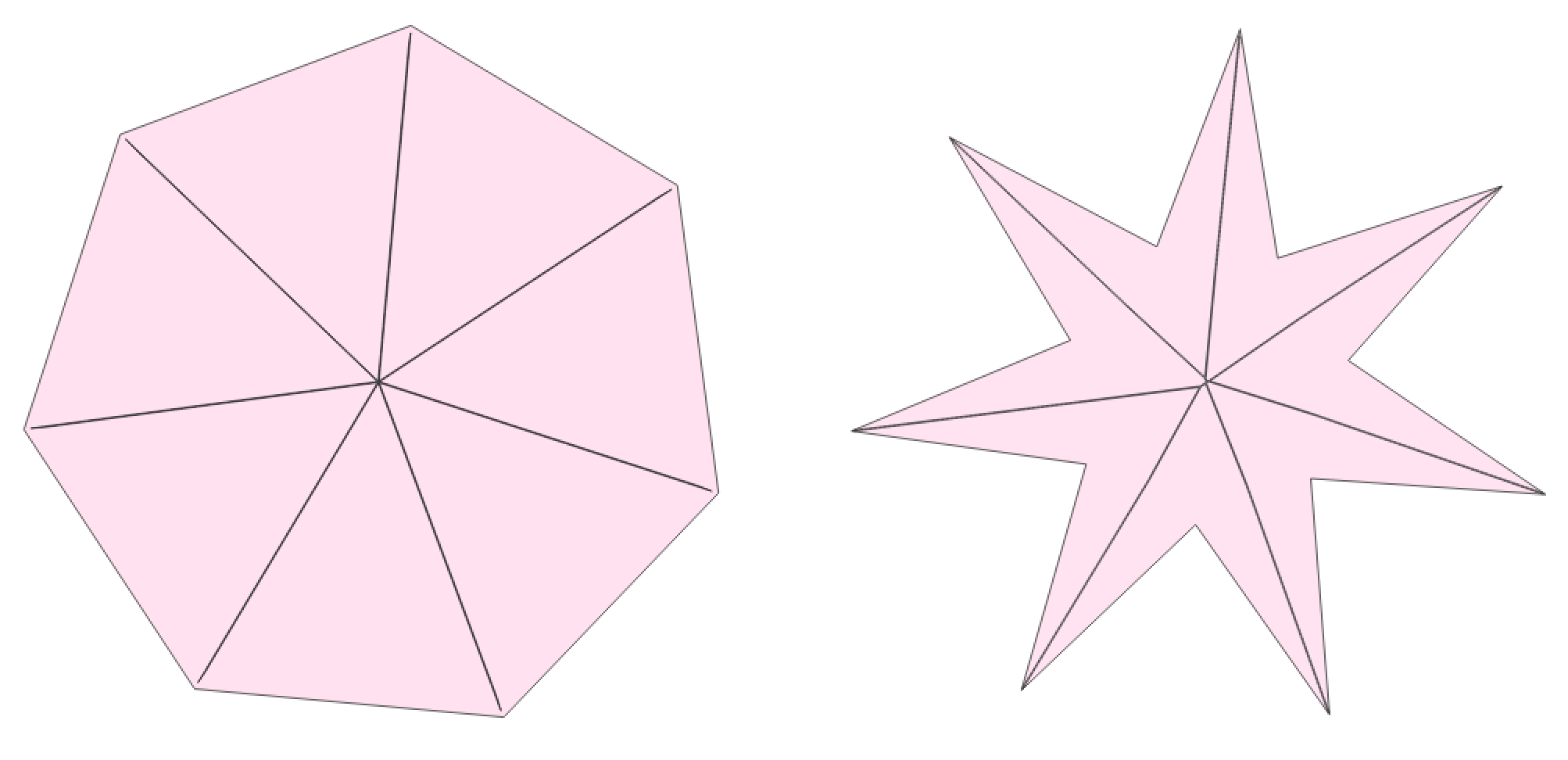}}
\caption{\label{fig:hep_sta} A heptagon and a 7-star (heptagram) with their medial axes: the astute reader will observe that these shapes share a medial axis, which is composed of line segments from the center which bisect convex angles.}
\end{figure}

We now detail some simple shapes and their medial axes. Honing an eye for where medial axes lie within shapes will be a useful tool as we move to analyze more complex figures. In Figures \ref{fig:squ_tri} through \ref{fig:hep_sta}, we detail the medial axes of a square, a triangle, a rectangle, a circle, a letter `E', a heptagon, and a 7-star.

\subsection{Voronoi diagrams}

Next, we will remind readers of the Voronoi diagram, a construct which eventually relates to the medial axis and will serve as the crux of our custom models for Congressional district analysis.

\begin{definition}
The {\sl Voronoi polygon\footnote{In the field of geographic information science, Voronoi polygons are generally referred to as Thiessen polygons. In this paper, we use Voronoi polygons.}} associated with an element $e$ from a set $S$ in the plane is the locus of points closer to $e$ than any other element of $S$. The collection of Voronoi polygons for each element of $S$ is called the {\sl Voronoi diagram} of the set $S$. 
\end{definition}

\begin{definition}
Given a simple polygon, a vertex $q$ is called {\sl convex} if the internal angle at $q$ is less than $\pi$. Otherwise, the vertex is called {\sl reflex}.
\end{definition}

We observe that the boundaries, called Voronoi edges, of Voronoi polygons are not always line segments. To see why some boundaries may be portions of parabolic curves, consider the construction of a parabola (the locus of points equidistant to a fixed point and a line). Further, we note that if a vertex touches one of the medial axis’ maximal circles, that vertex must be reflex.

The following string of theorems is extracted from \cite{planar}, which the reader may access for proofs or full development of this machinery. In each of these, $G$ is the planar shape whose medial axis we aim to construct.

\begin{lemma}
The Voronoi polygons $V(e_i)$ and $V(e_j)$ share an edge if and only if there exists a point $z$ such that the circle centered at $z$ with radius $d(z,e_i)=d(z,e_j)$ does not include any boundary point of $G$ in its interior\footnote{Each element $e_i$ and $e_j$ may be either a point or a line.}.
\end{lemma}

\begin{lemma}
Let $z$ be any point on a Voronoi edge $\bar{B}(e_i,e_j)$ of $V(e_i)$. The circle centered at $z$ with radius\footnote{$I(z,e_i)$ denotes the image of a point $z$ projected onto an element $e_i$, which is either a point or a line.} $d(z,I(z,e_i))$ is tangent to elements $e_i$ and $e_j$ at points $I(z,e_i)$ and $I(z,e_j)$.
\end{lemma}

\begin{theorem}
The medial axis of $G$ is the set of Voronoi edges\footnote{These edges are taken from the Voronoi diagram generated by the region's boundary lines {\sl and} vertices.} less the edges incident with reflex vertices.
\end{theorem}

To understand this intuitively, consider that Voronoi edges represent points equidistant to at least two different lines on a region's boundary. In this sense, it becomes clear that the medial axis is a subset of the Voronoi edges. Now consider the points which compose any Voronoi edge incident to a reflex vertex. These points lie equidistant to the lines incident to this vertex. The maximal circles centered at any of these points will clearly intersect the region's boundary at the reflex vertex and nowhere else locally, so they are in general not included the medial axis.

\begin{figure}
\centerline{\includegraphics[height=2in]{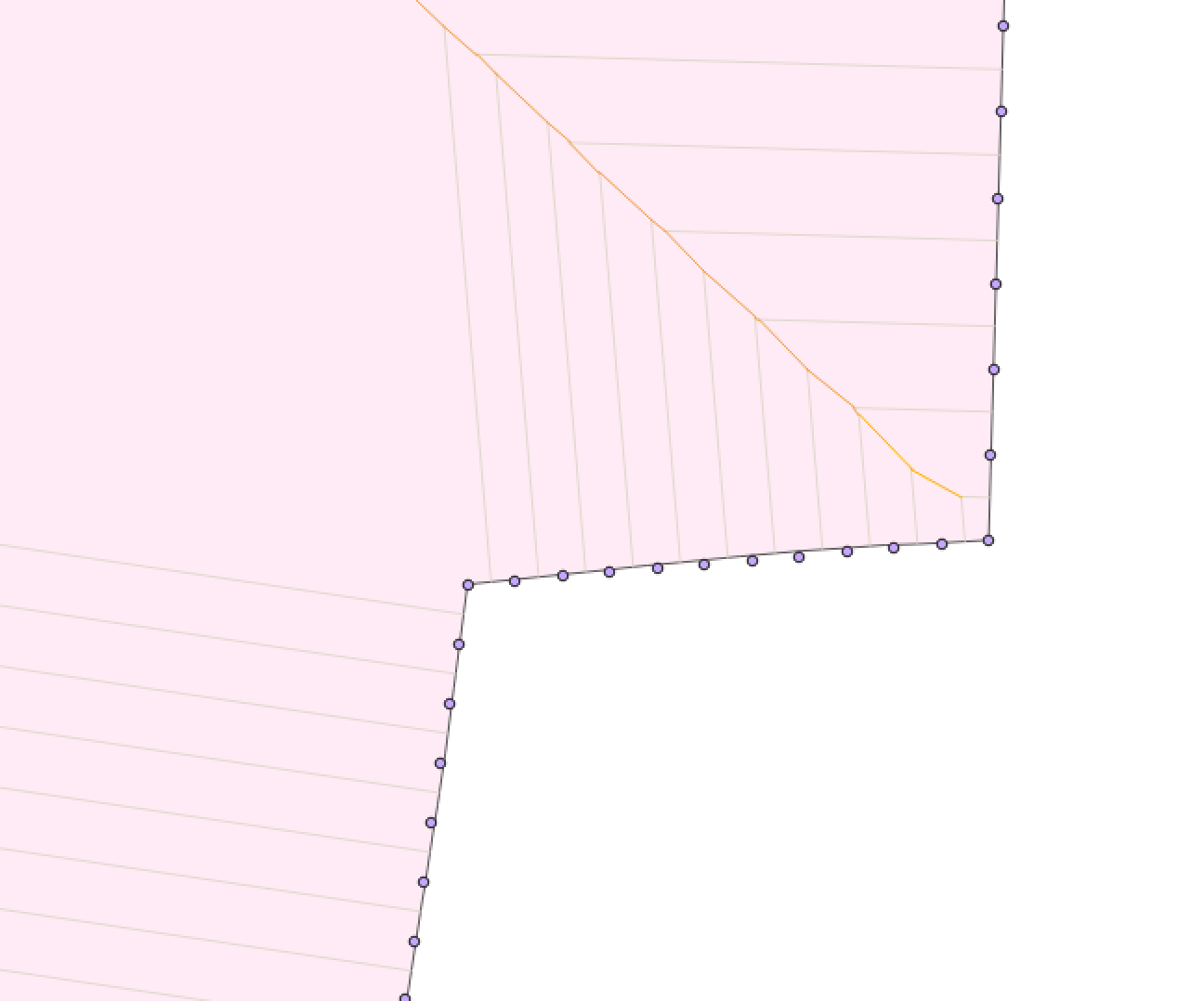}}
\caption{\label{fig:PA7d1} A detail of PA 7 (pink); orange lines compose our approximate medial axis and grey lines are those from our Voronoi diagram which intersect the boundary line and are thus deleted.}
\end{figure}

For computational practicality in this paper, we will approximate the Voronoi diagram of lines by taking the Voronoi diagram of discrete points populating these lines. As the number of points sampled from each line increases, a subset of our Voronoi edges approximates the Voronoi edges from the diagram taken on lines and thus approximates the medial axis. For diagrams taken on points, we observe that the Voronoi edges are line segments (i.e. none are portions of parabolic curves, since these result only from equidistance to a fixed point and a {\sl line}). Then, if we remove the Voronoi edges lying incident to a line on the region's boundary, we observe our approximation of the medial axis. An example of this from Pennsylvania district 7 is shown in Figure \ref{fig:PA7d1}.

\begin{remark}
Since our Voronoi diagrams are taken solely on points, our Voronoi edges will never intersect vertices---reflex or otherwise. We will only witness intersections to the boundary line.
\end{remark}

\begin{figure}
\centerline{\includegraphics[height=3in]{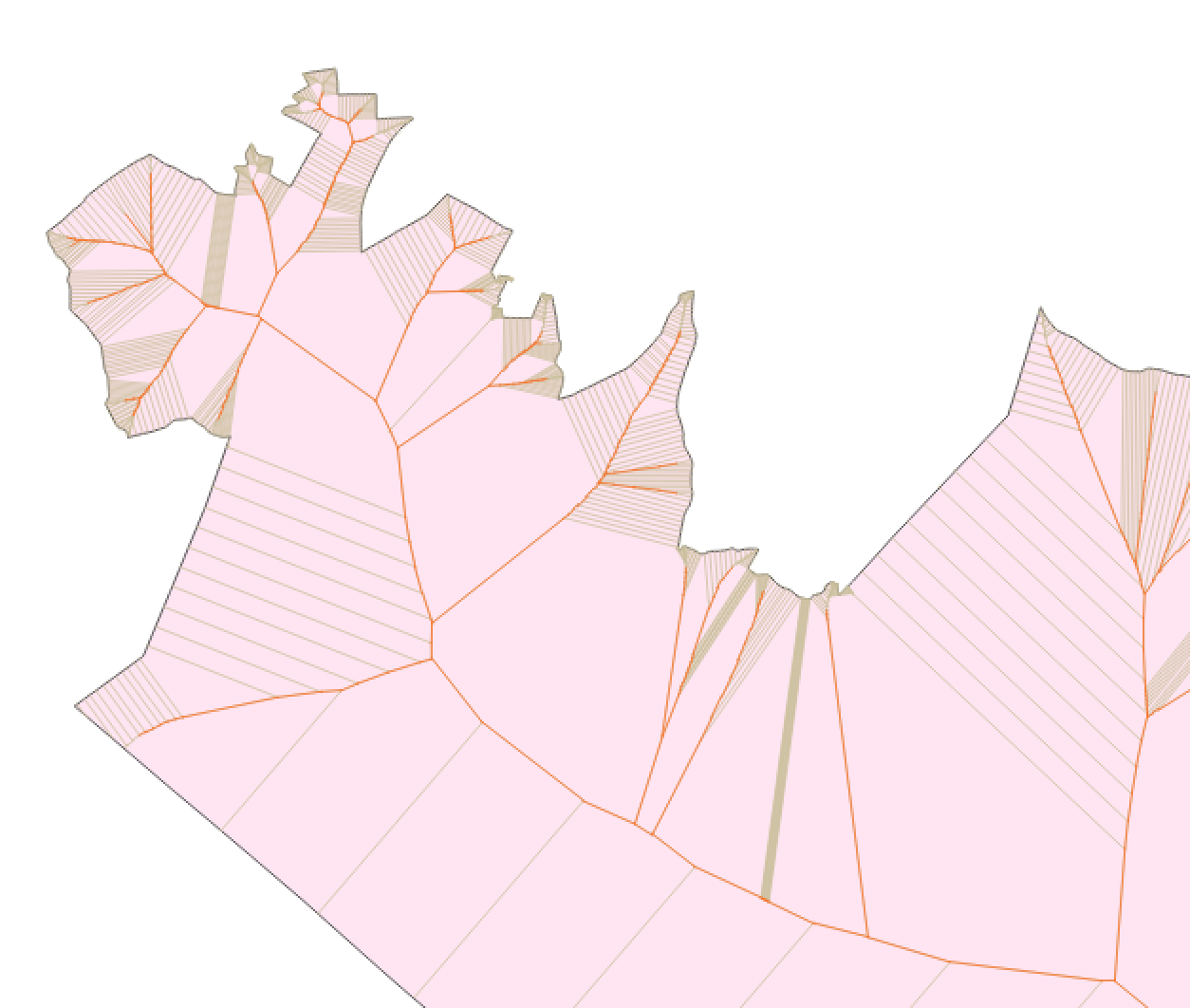}}
\caption{\label{fig:PA7d2} A detail of PA 7; orange lines compose our approximate medial axis and grey lines are those from our Voronoi diagram which touch the buffer zone and are thus deleted.}
\end{figure}

Thankfully, the margin of approximation becomes sufficiently small with even just a handful of points on each line. As described in Section \ref{sec:methods}, we will soon populate each Congressional district boundary line with 10 points, then take the Voronoi diagram of these points. By employing a buffer zone for extracting our medial axis from the Voronoi edges, we remove all edges touching the boundary (at lines, convex vertices, and reflex vertices). The Voronoi edges lying incident to lines and reflex vertices should indeed be removed to discern the medial axis. In this process, however, we also prune the tips of the medial axis where {\sl should} touch the boundary, as our buffer zone clips the outermost Voronoi edges (even those lying incident to convex vertices). This provides motivation for higher point densities---due to the line segmentation on the Voronoi diagram in our GIS software, with more points populating the boundary lines, we trim off less from the medial axis. Indeed, only the tip of this medial edge is removed---we are left with segments equidistant to the line segments second-closest to the convex vertex, and beyond. After pruning, our approximate medial axis is what remains. An example of this from Pennsylvania district 7 is shown in Figure \ref{fig:PA7d2}.

\begin{figure}
\centerline{\includegraphics[height=3in]{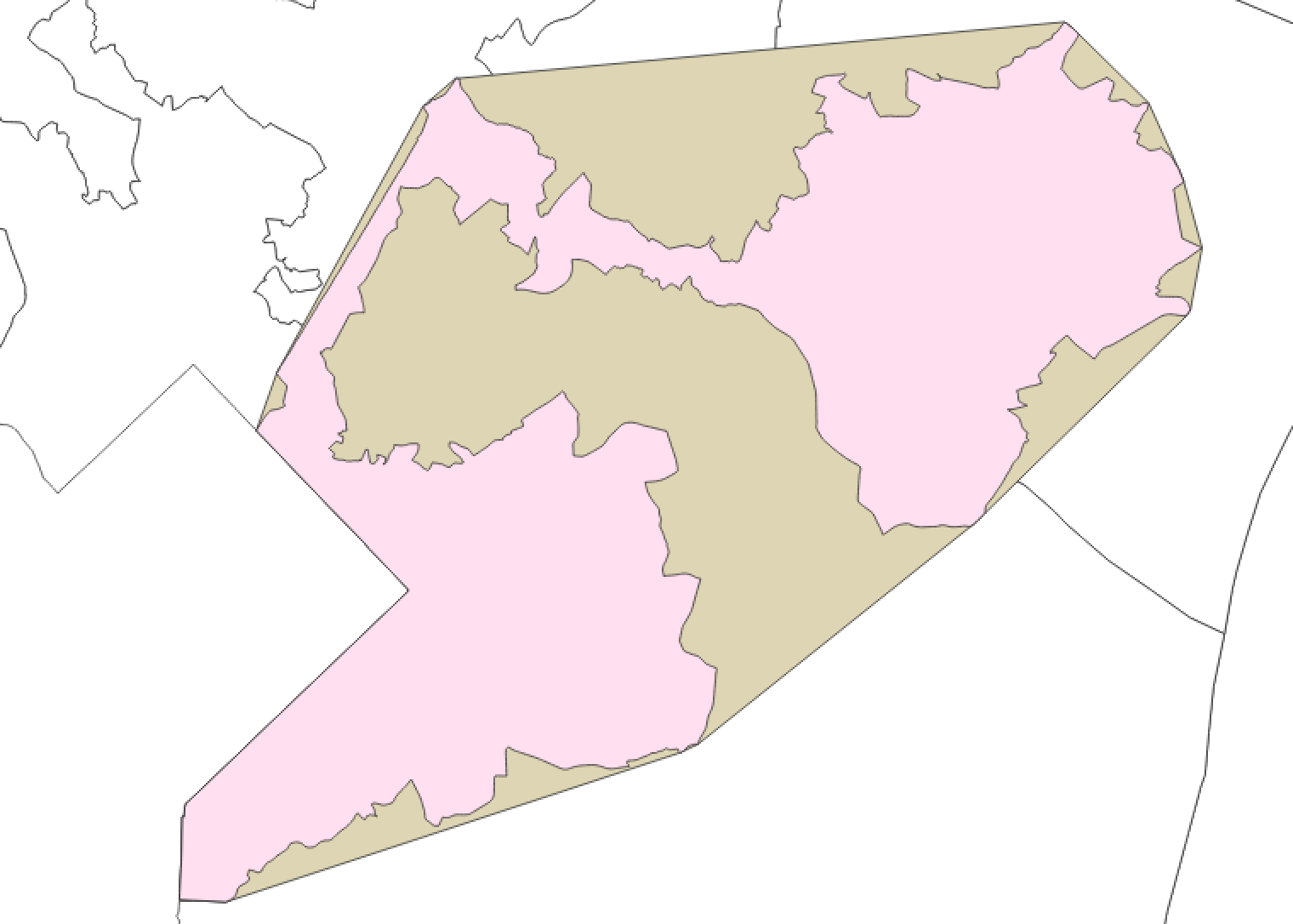}}
\caption{\label{fig:MD4h} Maryland district 4 (pink) displayed over its convex hull (tan) clipped to the Maryland state border.}
\end{figure}

Lastly, we remind readers of another definition, which we will soon need in our construction of the medial-hull ratio. This concept is demonstrated in Figure \ref{fig:MD4h}.

\begin{definition}
The {\sl convex hull} of a region $R$ is the intersection of all convex sets containing $S$. For our purposes, it is the smallest convex polygon containing $R$.
\end{definition}

\section{Methods}\label{sec:methods}
In order to compute the medial axes of Congressional districts, we employ QGIS 3.0 Girona.

Since Congressional district boundaries form simple polygons, the true medial axis of each district is the set of border lines, from Voronoi polygons generated on boundary lines, which do not intersect reflex vertices. We remark again that what we compute and discuss in this paper is not precisely the medial axis, but rather an approximation by way of a Voronoi diagram, taken on vertices populating boundary lines and then clipped and pruned. As the number of points populating each boundary line approaches infinity, our approximation grows closer to the true medial axis; however, computation time increases outrageously. In this paper, we refer to our approximation of the medial axis as simply ``our medial axis.''

In this section, we detail the procedure for handling and processing Congressional district data. Once the data is collected and prepared for our models, we will compute the medial axes both of districts and of their convex hulls clipped to state borders. We then will compute the total length of these axes and divide them to procure the medial-hull ratio, which is finally our tool for normalized comparison between districts and states.

\subsection{Data}
Almost all computations in this paper are made on data provided by The United States Census Bureau as TIGER/Line Shapefiles for the 115\textsuperscript{th} Congressional Districts \cite{TIGER}. This dataset is considered the most comprehensive dataset available and was made public on September 28, 2017. The other data source used is the 2018 remedial Congressional map of Pennsylvania\footnote{The remedial plan was adopted on February 19, 2018 as a result of the case {\sl League of Women Voters v. Commonwealth of Pennsylvania} (2018), which claimed that since ``a diluted vote is not an equal vote,'' the old district map violated the state Constitution \cite{newpa}.}. The remedial Shapefile is available via The Unified Judicial System of Pennsylvania \cite{newpa}.

\subsection{Reprojection}
In order to make meaningful distance and area calculations on Shapefile data, we use GIS software to project the data into a proper coordinate reference system (CRS). For our purposes, reprojecting data is usually handled by state, but may be applied to a collection of districts best covered by the same Universal Transverse Mercator (UTM) zone. For example, the western and eastern portions of Pennsylvania traverse the 17N and 18N UTM zones, respectively. We note that calculations are reasonably accurate in regions just outside the proper UTM zone, so we occasionally treat states uniformly for convenience. An instance of this is Illinois, which falls mostly in the 16N UTM zone, though its western edge crosses over to the 15N UTM zone.

Once the desired districts and their relevant UTM zone are decided, the region should be selected in QGIS from its originating dataset (here the ``Select Features By Value'' tool is useful if an entire state is sought in isolation). The sourced data is projected in NAD83 (EPSG:4269), a North American geodetic datum. The user should copy the selected features, then paste them as a new vector layer. When prompted to save this new layer, select the desired CRS (the present authors used UTM zones, but there are others for which distance calculations have meaningful dimension). After saving the new layer as its own Shapefile, close QGIS and restart the program by opening the saved Shapefile. The desired layer and project-wide CRS should now render correctly with coordinate units now displayed (in our case, meters).

\subsection{Processing}
The state vector layer, now correctly projected, should be separated (via {\sl Split vector layer} with the 115th Congressional district FIPS code ``CD115FP'' as the unique ID field) so that each district in the region corresponds to a unique layer. For this paper, we developed a custom {\sl Compute medial axis} model to take any district vector layer and output its medial axis as a new vector layer.

\subsubsection{Custom model: Compute medial axis}

\begin{figure}
\centerline{\includegraphics[height=3in]{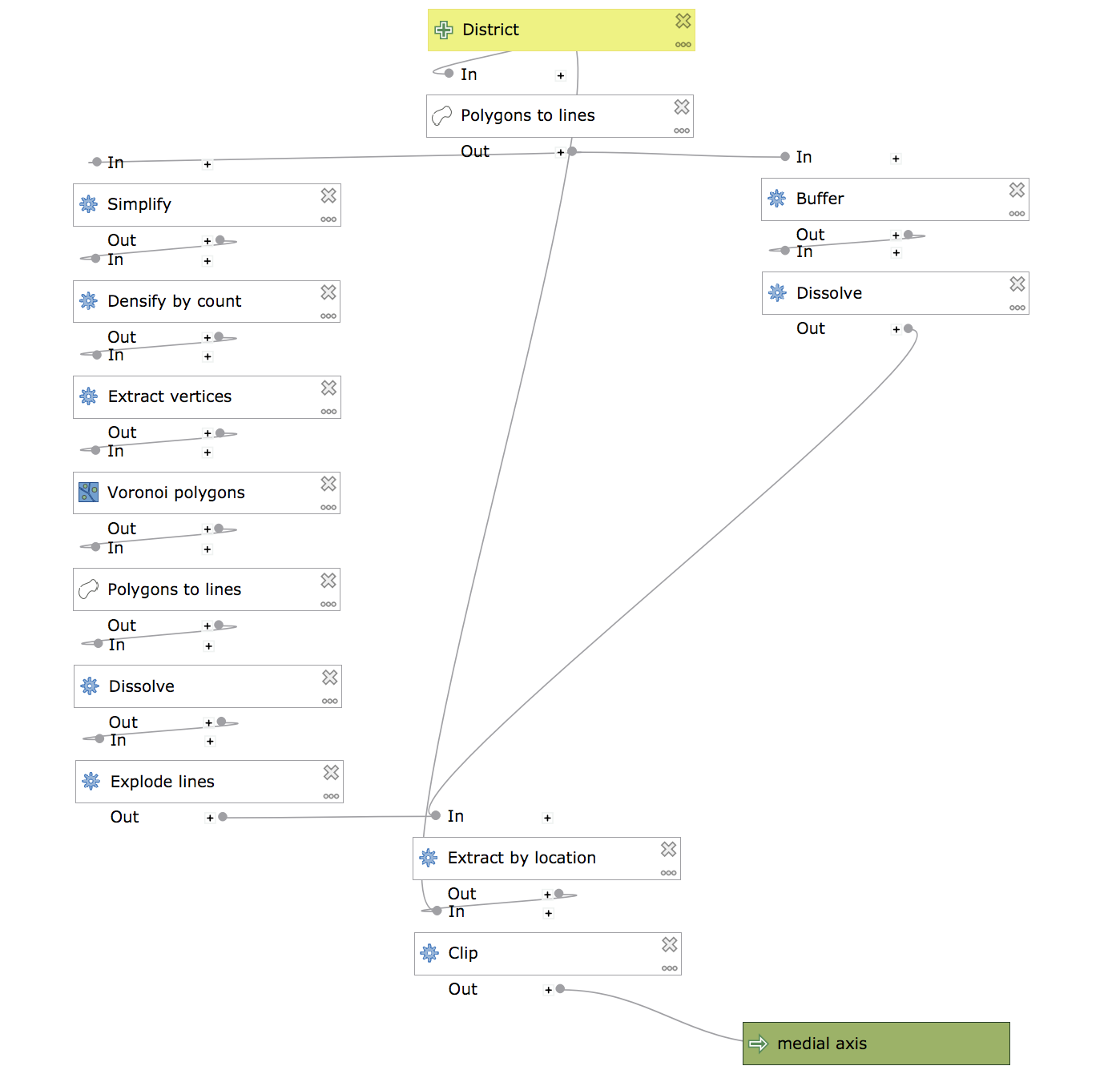}}
\caption{\label{fig:CMA} Diagram of our custom {\sl Compute medial axis} model in QGIS; the input is a `District' polygon layer and the output is a `medial axis' line layer.}
\end{figure}

With a district input as a vector layer (with polygon geometry), the model applies {\sl Polygons to lines} to extract the district boundary. Then, with {\sl Buffer} (distance parameter\footnote{For our custom models, assume that parameter values are the QGIS 3.0 default values, unless otherwise stated.} set to 200). The parameter value was chosen empirically and corresponds to a buffered zone extending 200 meters inward (and also outward) from the district boundary. Imperfect alignment of disparate QGIS geometries demands that a buffer zone be used, and the present authors found that this parameter value forgives sufficiently small boundary noise while accommodating for the fact that the forthcoming Voronoi diagram will be computed on points rather than lines (as in the true medial axis; see Section \ref{sec:medial}). To satisfy pipeline restrictions in the QGIS Processing Modeler, our model then applies {\sl Dissolve} to the buffer zone so that its features are consolidated into a single entity.

Meanwhile, the model manipulates the district boundary via {\sl Simplify} (tolerance parameter set to 500 for the Douglas-Peucker algorithm\footnote{This algorithm was first described in \cite{DP}.}). This parameter value was also chosen empirically and judged to grant an acceptable trade-off between computational ease and district shape integrity. By negligibly simplifying the shape, the same small boundary noise is forgiven as with the buffer zone, and the Voronoi diagram is then taken on much fewer points to allow for considerably speedier processing times.

Compromising now with this simplification, the model then densifies the simplified boundary line with {\sl Densify by count} (vertex parameter set to 10). Increasing the parameter value here allows the Voronoi diagram on these points to more closely approximate the Voronoi diagram taken on boundary lines (which would contain the true medial axis as a subset). However, computation quickly becomes cumbersome, as discussed above. This value was chosen as a moderate compromise between these conflicting motivations. Minor issues remain only near relatively long boundary lines, which are the most desperate for these densified points since the most ostensible discretion between a line and finite points lies here. This spacing issue would be solved by implementing the {\sl Densify by interval} algorithm from the QGIS Processing Toolbox, but the present authors witnessed pipeline errors stemming from this alternative.

The model then prepares with {\sl Extract vertices} from the simplified-then-densified boundary before computing the {\sl Voronoi polygons} of these vertices. Ideally in this step, the Voronoi polygons would be computed for the boundary lines, as discussed above, but due to computational and software restraints, we compromise with vertices. The Voronoi diagram created extends to the coordinate-based extent of its generating vertices. Our model applies {\sl Polygons to lines} now to the Voronoi diagram to yield the borders of the Voronoi polygons, then employs {\sl Dissolve} and {\sl Explode lines} successively to these outputs in order to account for the doubling up of lines bordering two different Voronoi polygons. Our medial axis is now close by.

In the final moments, the model recalls the dissolved buffer zone from initial stages in order to filter out Voronoi lines which are not a part of our medial axis. As discussed above, this procedure removes Voronoi lines which stem from reflex vertices or boundary lines. Unfortunately, it also pares off a small portion of the medial lines which touch convex vertices. Users will witness this negligibility, which will anyway be divided out with the eventual ratio metric actually employed for analysis. The model uses {\sl Extract by location} (feature parameter set to ``disjoint'' on the exploded lines against the dissolved buffer) then prunes out lines which lie beyond the district boundary with {\sl Clip}\footnote{The {\sl Clip} algorithm can be utilized as early as immediately following {\sl Voronoi polygons}. This was noticed to be the ideal location for this algorithm in the pipeline with respect to computation time.}.

Now, we witness our medial axis, which is saved as a new vector layer by our model. Computation of the model lasts roughly 10 seconds for a single district via QGIS 3.0.1 on a 2016 MacBook possessing an Intel Core m7 processor with 1.3 GHz speed and 8 GB of memory.

\subsubsection{Computing the hull axis}
In addition to computing the medial axis of the district itself, we then compute the medial axis of the district's convex hull after having been clipped to the district's state border. Motivation for this will be discussed below. This process is similar, save for a few additional preparation steps before we invoke our custom {\sl Compute medial axis} model once again.

We accept an additional input geometry here: a vector layer containing every district in the state to which the district in question belongs. Once this layer is prepared, we derive the {\sl Convex hull} of the district, then {\sl Clip} the hull against the state input layer. Motivation for this lies in avoiding penalty for erratic or otherwise unwieldy state boundary lines. Manipulation of these borders lie beyond the powers of Congressional district cartographers---those who may stand accused of gerrymandering. Then, after clipping, we pipe the hull through {\sl Compute medial axis} and extract the medial axis of the district's clipped convex hull, which we may refer to simply as the district's hull axis.

When this hull axis procedure is formalized into another custom model, computation lasts roughly 8 seconds for a single district under the same circumstances as described above. While we embark on additional steps before the {\sl Compute medial axis} model, these steps lead to a much simpler input to that model, thus reducing overall computation time.

\subsubsection{Computing the medial-hull ratio}

Let $\mathcal{M}(P)$ denote the medial axis of a polygon $P$ and $\mathcal{C}(P)$ denote its clipped convex hull. Let $|A|$ denote the total length of an axis $A$.

\begin{definition}
The {\sl medial-hull ratio} of a polygon $P$ is defined to be $$\mathcal{R}_P = \frac{|\mathcal{M}(P)|}{|\mathcal{M}(\mathcal{C}(P))|}.$$
\end{definition}

Building upon our previously detailed models, we wrote a consolidated {\sl Compute ratio} model to streamline computation and length summation for these medial and hull axes. The procedure for calculation of this medial-hull ratio is detailed below.

The inputs, as with the hull axis model, are two vector layers: one containing the district of interest and the other containing that district's entire state. Taking the district, the model computes its medial axis as above. With the output, it then calculates the distance over which the medial axis stretches, using {\sl Sum line lengths} over the district. To placate QGIS pipeline restrictions again, the model uses {\sl Dissolve}. The calculated total length is then added as an attribute to the original district by way of {\sl Refactor fields}. In this step, all attributes are deleted except for the state FIPS code (``STATEFP''), the 115th Congressional district FIPS code (``CD115FP''). We append an additional field (``Medial'') and populate it with the output (``LENGTH'') from {\sl Sum Line Lengths}. We note that the refactored field now bears the shape of the originating district, which is desirable, and is attributed by only the three fields listed above.

Next, the model treats the refactored district as an original, and uses the input state layer to call our other model to compute the medial axis of the district's convex hull, as detailed above. In similar fashion as before, the model then uses {\sl Sum line lengths} before {\sl Dissolve} and {\sl Refactor fields}. This time, we append two additional fields (``Hull'' and ``Ratio,'' respectively) and populate them with our most recent output (``LENGTH'') from {\sl Sum Line Lengths} and with the medial-hull ratio calculation (``Medial'' / ``LENGTH''), respectively. Finally, the model clips the output, which bears shape of the clipped convex hull, to the original district. In the end, we are left with a vector layer bearing shape of the desired district and possessing five attributes: ``STATEFP,'' ``CD115FP,'' ``Medial,'' ``Hull,'' and ``Ratio.''

\begin{remark}
In QGIS, the ``Medial'' and ``Hull'' fields are displayed in meters, as dictated by the CRS we used. The medial-hull ratio (``Ratio''), then, is a dimensionless measure that may be utilized for comparison of districts both within a given state and across state borders.
\end{remark}

This medial-hull ratio value is what is considered for analysis in this paper, both with respect to individual districts and to statewide averages. The medial-hull ratio is constructed by comparing a district's medial axis to that of its clipped convex hull because the convex hull, barring population restrictions, is one option of what cartographers could have used to generally compactify a district within the extent of land already covered.

\subsection{Issues}

We remind the readers of three choices that determine this incarnation of our model: the buffer parameter, the Douglas-Peucker algorithm parameter, and the densify count parameter. A model that allows these parameters to vary according to some measure of a district's size would largely eliminate size bias.

Moreover, future work should consider a measure on the count of connected components in each medial axis. In general, the medial axis of a connected component is itself connected, but due to our implementation with a buffer zone, regions where the district becomes thinner than 400 meters divides our medial axis. So, this would be a measure on frequency of thinness. Stronger yet would be other measures of graph complexity of the medial axis as a tree (or, indeed, a forest). Another avenue might be to compare the length of the smallest non-trivial `cut point' to the overall diameter of a district.

\subsubsection{Combs}

\begin{figure}
\centerline{\includegraphics[height=2.0in]{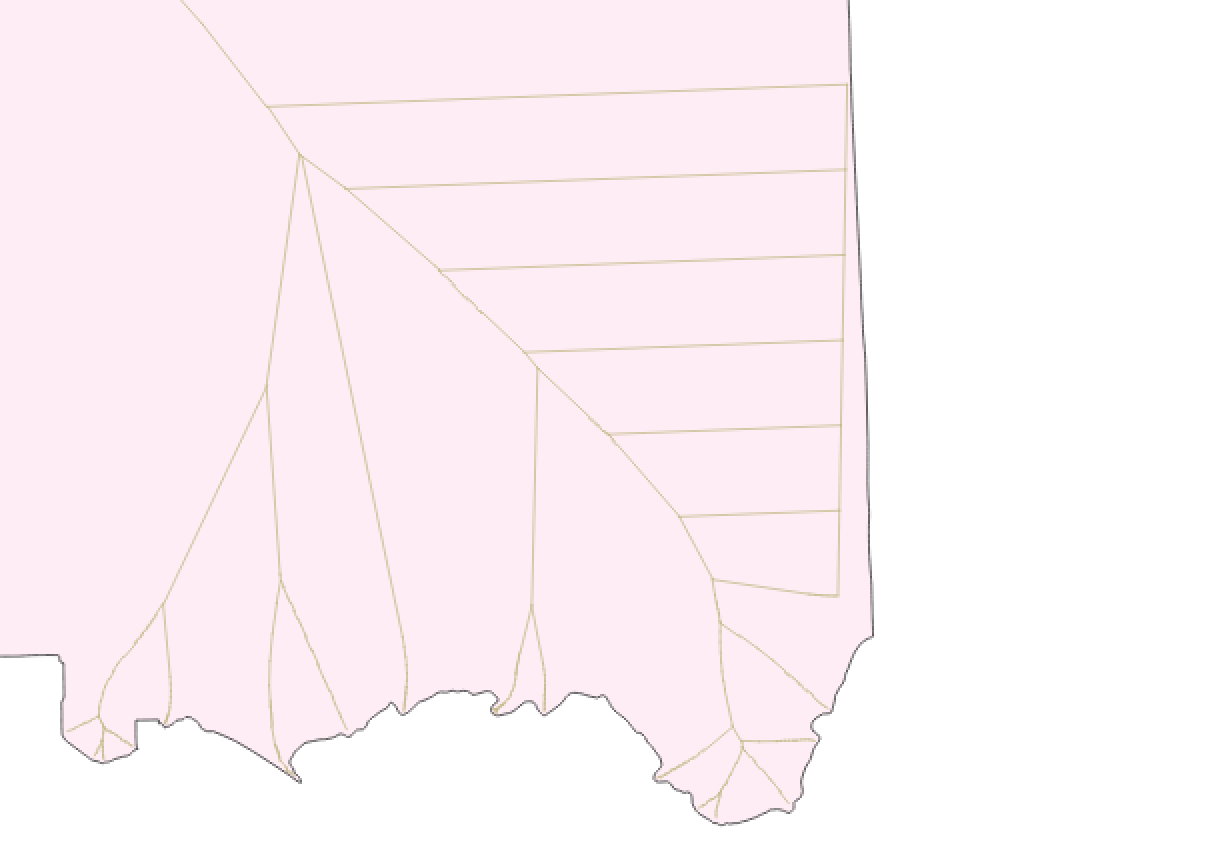}}
\caption{\label{fig:NV04d} Detail of comb in NV 04}
\end{figure}

In districts such as Nevada district 4 (Figure \ref{fig:NV04d}) and Florida district 27, the present authors have witnessed errors in the computed Voronoi polygons within the {\sl Compute medial axis} model. These resemble small combs, and carry through to the final output of medial axis. We suspect this is due to the way in which QGIS clips a set of Voronoi polygons to the rectangular extent of its generating vertices. When borders are simplified in our model, this may slightly shrink the extent, while the buffer zone surrounds the true border. Then, lines on the simplified border trespass on the medial axis by not being sheared off by the buffer. These vestiges seem small, seldom, and only marginally disruptive to the total medial axis length computed.

A possible solution might be to add distant vertices before computing the Voronoi diagram, thus widening the vertex extent and hopefully pushing any combs outside of the district.

\subsubsection{Infinite fragments}
In districts such as Illinois district 8,  the present authors have witnessed errors of unknown origin. These are geographically small fragments in the medial axis which possess near-infinite length by calculation. In this sense, the {\sl Sum line lengths} algorithm is disrupted and returns meaningless output for the medial axis' total length, though the issue is very insular.

No solutions have been posed, though a workaround exists by manual excision after trial-and-error location of the fragment. What remains should bear computationally sound length.

\subsubsection{Longer hull axes}\label{subsec:longHull}
In districts such as WV 3, the present authors have witnessed an unintuitive consequence of our processing paradigm. Though the hull axes are in general simpler than medial axes, sometimes their length will be greater. This is best understood when conceptualizing the axis with skeletal terminology---the medial axis possesses limbs which poke out from some deeper spine. Consider the case in which a district lies on a jagged state border. The medial axis will bear many limbs stretching out to this jaggedness, as expected. However, where the hull axis normally would smooth over this jagged side of the district, it cannot do so since the jaggedness is derived from the state border. Yet, the opposite side of the district may be smoothed over and simplified by the convex hull, moving the spine further from the jagged side. Then, the hull axis possesses a similar count of limbs, but each of which is now longer as it must stretch further out toward that jagged state border. In cases (especially of near-convexity) where this added length surpasses the saved length from axis simplification, the medial-hull ratio is then less than 1.

A possible solution might be to eschew penalty for immutable state borders in some other way. This might include smoothing over (in style of the convex hull) every district boundary which lies on a state border, even for the medial axis computation of a district. This would avoid penalizing jagged district boundaries which stem from state borders while at first glance avoiding the current issue\footnote{The fact that certain districts bear medial-hull ratio less than 1 should not cause duress. While evidencing gerrymandering should be done via high medial-hull ratios, the present authors make no claim that low ratios indicate a strong lack of political bias}.

\subsubsection{Sensitivity to noise}\label{subsec:noise}
Perhaps the deeper issue from that listed above is sensitivity to noise. Small juts, especially those in otherwise smooth boundary regions, add considerable length. Depending on how `thick' the district is near that region, even just one jut could earn a heavy penalty to the medial axis length. In certain cases, this may penalize the medial-hull ratio too harshly for how much malicious political intent might have motivated the carving of that jut. Indeed, the jut must be at least 400 meters in width to filter past the buffer zone; changing the buffer parameter value in our model would affect how extreme these juts need to be in order to warrant penalty. Indeed, this leads to bias against districts with larger area. The penalty on their medial axes for small, distant boundary noise is high. Since this is smoothed over for the hull axis, the medial-hull ratios of large districts are perhaps inflated.

Implementing the buffer zone and simplifying the district boundary already work to forgive sufficiently small noise. A choice must be made on how small of noise to penalize and use as witness to gerrymandering. However, other solutions that penalize small, distant juts in the boundary less harshly should be sought in future work; in order to accurately flag evidence of gerrymandering, our focus should be balanced between the district's overall shape and the extremity of noise on its boundary.

\subsubsection{Convexity}
If a district itself is (near) convex except for state borders, its medial-hull ratio will be equal to 1. This district may still be gerrymandered, though it would slip through the cracks of analysis in this paper. Pragmatically, most convex districts would be well considered compact. A prime counterexample, however, would be to consider a long but skinny rectangular district which traverses an entire state: its medial-hull ratio would be 1.00, but this district would be very probably gerrymandered. No solution to this within medial analysis has been posed.

\subsubsection{Partial gerrymandering}
In such cases as PA 17, we witness circumstances where part of a district has a convex squarish shape while another portion has an jagged border which is almost surely the result of gerrymandering. However, the ultimate affect on the medial-hull ratio might not alert us to presence of this gerrymander: the effect of the `bad' region is partly masked by that of the `good.'

Possible solutions lie well beyond the scope of this paper, but may include sampling metrics akin to the medial-hull ratio for many subregions of a district, then weighting them to compose an ultimate metric for the entire district.

\section{Computational Experiments}\label{sec:computation}

In this section, we demonstrate specific results on a few states from the model described above. For Pennsylvania, we will go into the most detail, investigating each district's shape for the 115th Congressional map before addressing the improvements heralded by the remedial map. Next, we will analyze sample districts from North Carolina and Florida. After these states are detailed, we discuss comparisons made between states and how that may also evidence gerrymandering, which generally occurs on a state-wide basis instead of solely by district.

Here, we also propose a four-tier categorization system for medial-hull ratios. For a district D, the medial-hull ratio $\mathcal{R}_D$ is judged as follows:
$$\text{If }\mathcal{R}_D < 2.00\text{, then $D$ is Category 1.}$$
$$\text{If }2.00 \leq \mathcal{R}_D < 2.40\text{, then $D$ is Category 2.}$$
$$\text{If }2.40 \leq \mathcal{R}_D < 2.80\text{, then $D$ is Category 3.}$$
$$\text{If }2.80 \leq \mathcal{R}_D\text{, then $D$ is Category 4.}$$

For analysis, we consider Category 1 to bear no explicit evidence of gerrymandering\footnote{Caveat: we claim that districts in Category 1 have no explicit evidence of gerrymandering {\sl under our method of analysis}. We do not claim that this necessarily implies a lack of gerrymandering in these districts.}, Category 2 to bear suspect evidence of gerrymandering, Category 3 to bear probable evidence of gerrymandering, and Category 4 to bear near-certain evidence of gerrymandering.

Throughout this section, we may refer to districts by abbreviation of their state name followed by district number. For example, Pennsylvania district 5 would be denoted as PA 5.

\subsection{Pennsylvania}

Our first in-depth look details Pennsylvania\footnote{State FIPS code 42}. We acknowledge the importance of this state in the current legal scene surrounding gerrymandering: its 115th Congressional district map, which had been in place since 2011, was redrawn by the Supreme Court of Pennsylvania in early 2018. Here, we will investigate both the original and remedial maps.

\begin{remark}
We note that the eastern portion\footnote{Here, we refer to districts which generally fall within boundary of the 18N UTM zone: Pennsylvania districts 1, 2, 4, 6, 7, 8, 10, 11, 13, 15, 16, and 17.} of Pennsylvania served as the sandbox of sorts for the development of our custom {\sl Compute medial axis} model in QGIS. Specific focus was laid on Pennsylvania district 7 (PA 7), which inspired most of the parameter value choices detailed in Section \ref{sec:methods}. This district was chosen due to its modestly errant shape and somewhat average\footnote{This was computed by way of geometric mean on district areas of a few sample states.} geographic size.
\end{remark}

\begin{figure}
\centerline{\includegraphics[height=2.5in]{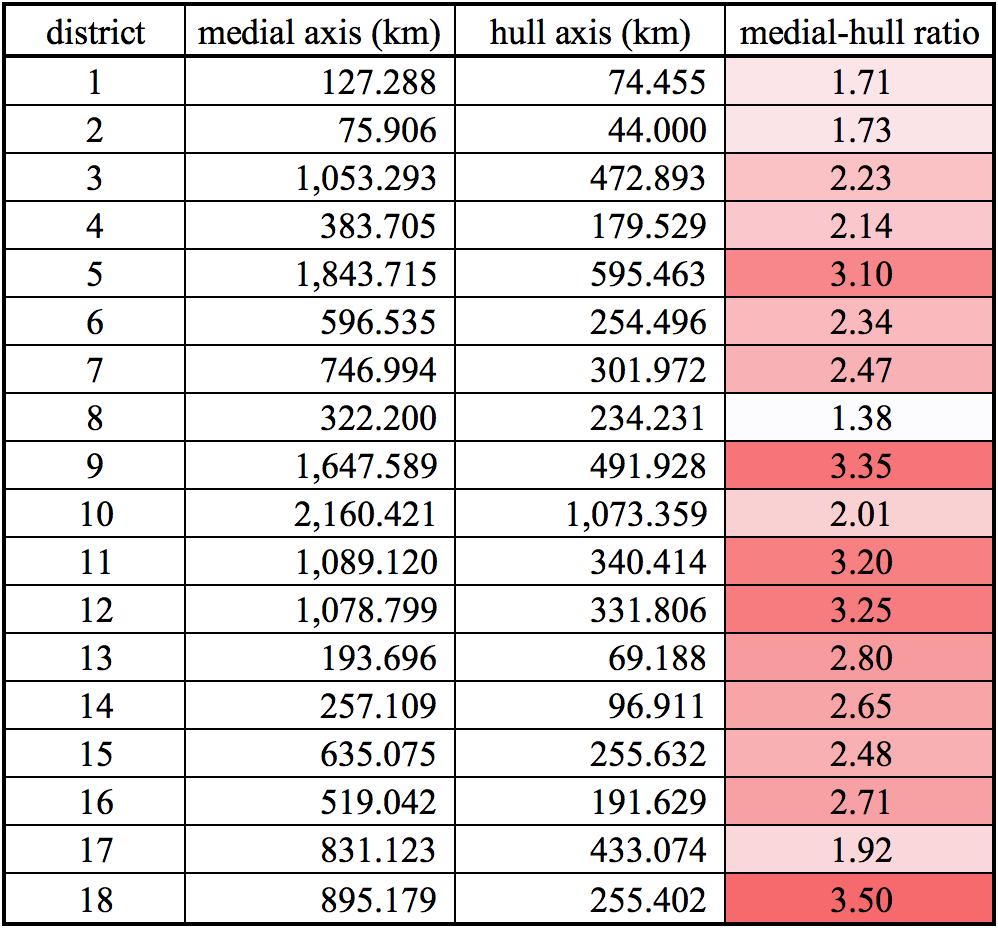}}
\caption{\label{fig:PA} Pennsylvania district data}
\end{figure}

\subsubsection{Pennsylvania district 1}

\begin{figure}
\includegraphics[width=0.475\textwidth]{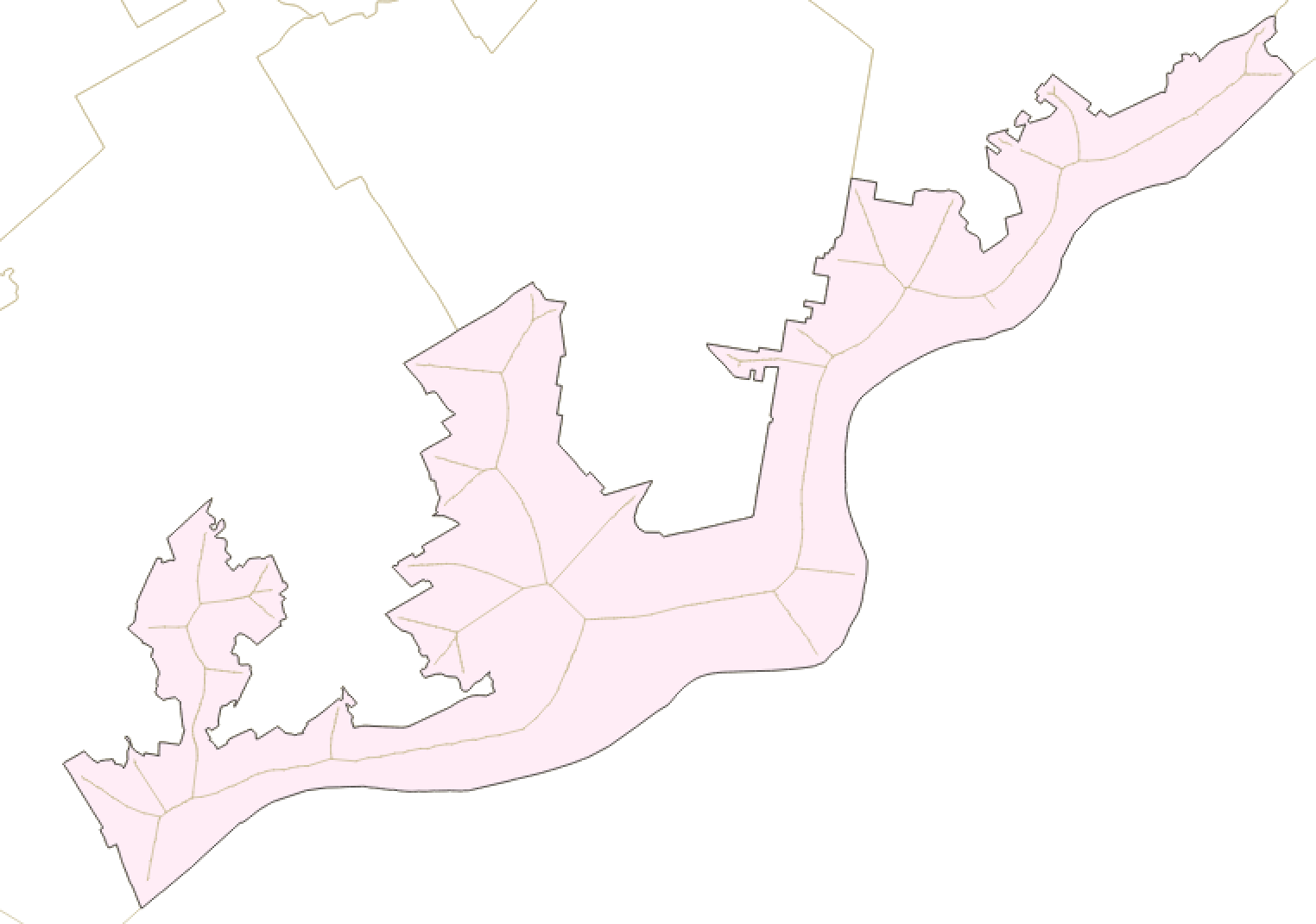}
\hspace{\fill}
\includegraphics[width=0.475\textwidth]{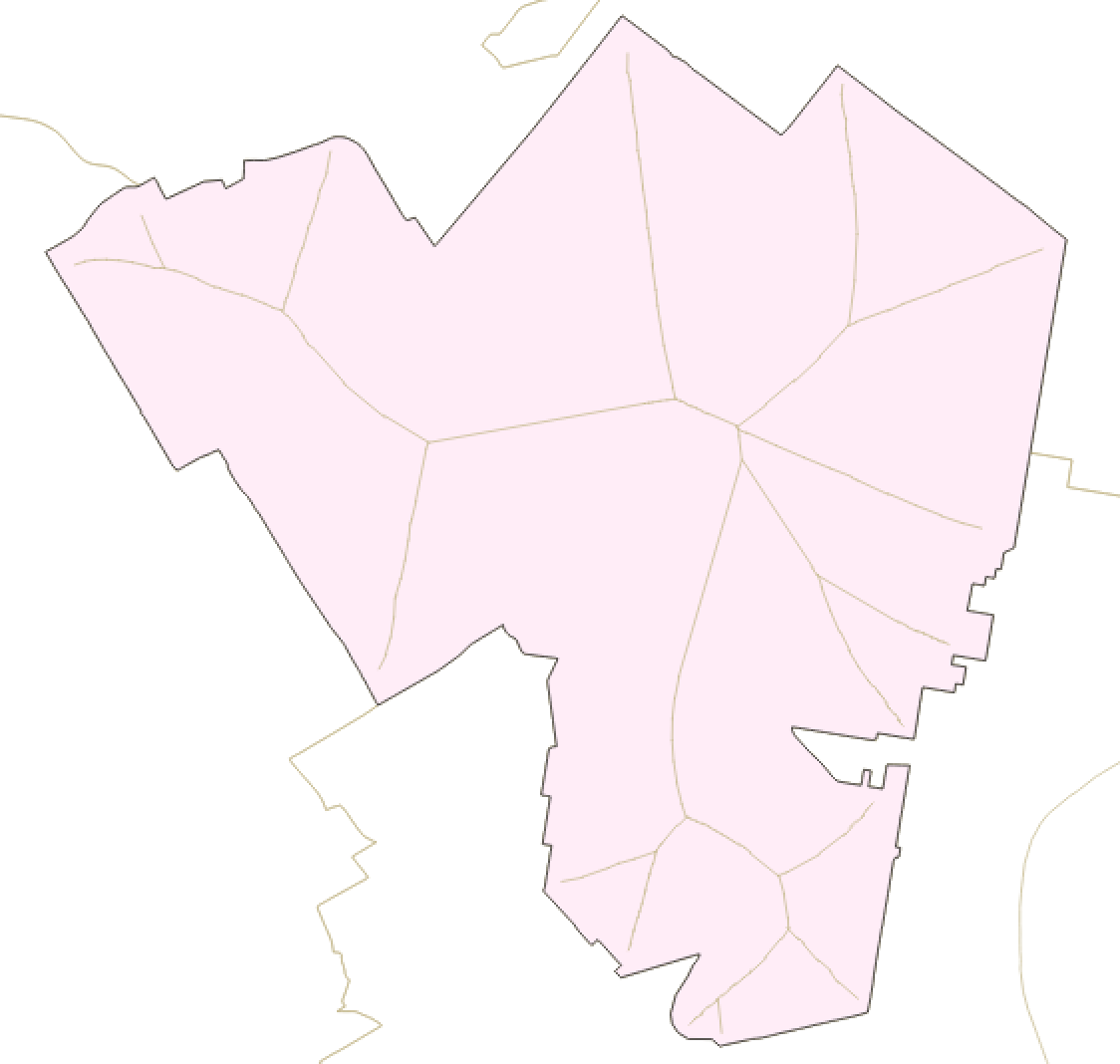}
\caption{Medial axes of PA 1 (left) and PA 2 (right)}\label{fig:PA01-02}
\end{figure}

This district, lying on a small portion of the state's southeast border, bears a generally thin rectangular shape that stretches across portions of Philadelphia as well as swaths of land both to its northeast and southwest. There is a small bubble erupting from the southwest portion of the district, leaving the area to its east excised from the district boundary. This analysis might fall victim to the issue detailed in Section \ref{subsec:longHull}, since the district has skinny sections along the state border (here subject to the winding Delaware river) which do not greatly increase the medial length though they seem to evidence gerrymandering. Indeed, this district possesses a medial-hull ratio of 1.71 (cat.~1), which is the second-lowest in its state. We remind the readers here that a lower ratio is not intended to claim there is a lack of gerrymandering seen; this metric only suggests where gerrymandering does exist, by way of high values.

\subsubsection{Pennsylvania district 2}
Here we witness a more modestly drawn district which is somewhat triangular and seemingly compact, save a few excisions on the southern end. The land covered includes West Philadelphia and a region north-northwest of the city. Indeed, the medial-hull ratio of 1.73 (cat.~1) is the third-lowest in Pennsylvania, and we witness no clear suggestion of gerrymandering from the district shape.

\subsubsection{Pennsylvania district 3}

\begin{figure}
\includegraphics[width=0.3\textwidth]{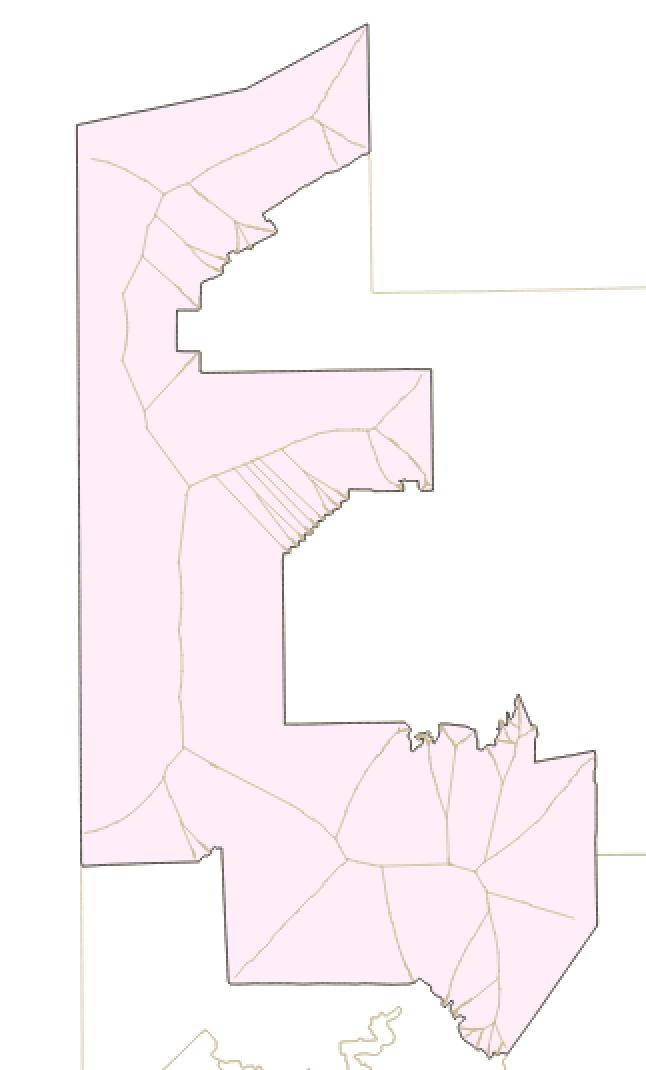}
\hspace{\fill}
\includegraphics[width=0.6\textwidth]{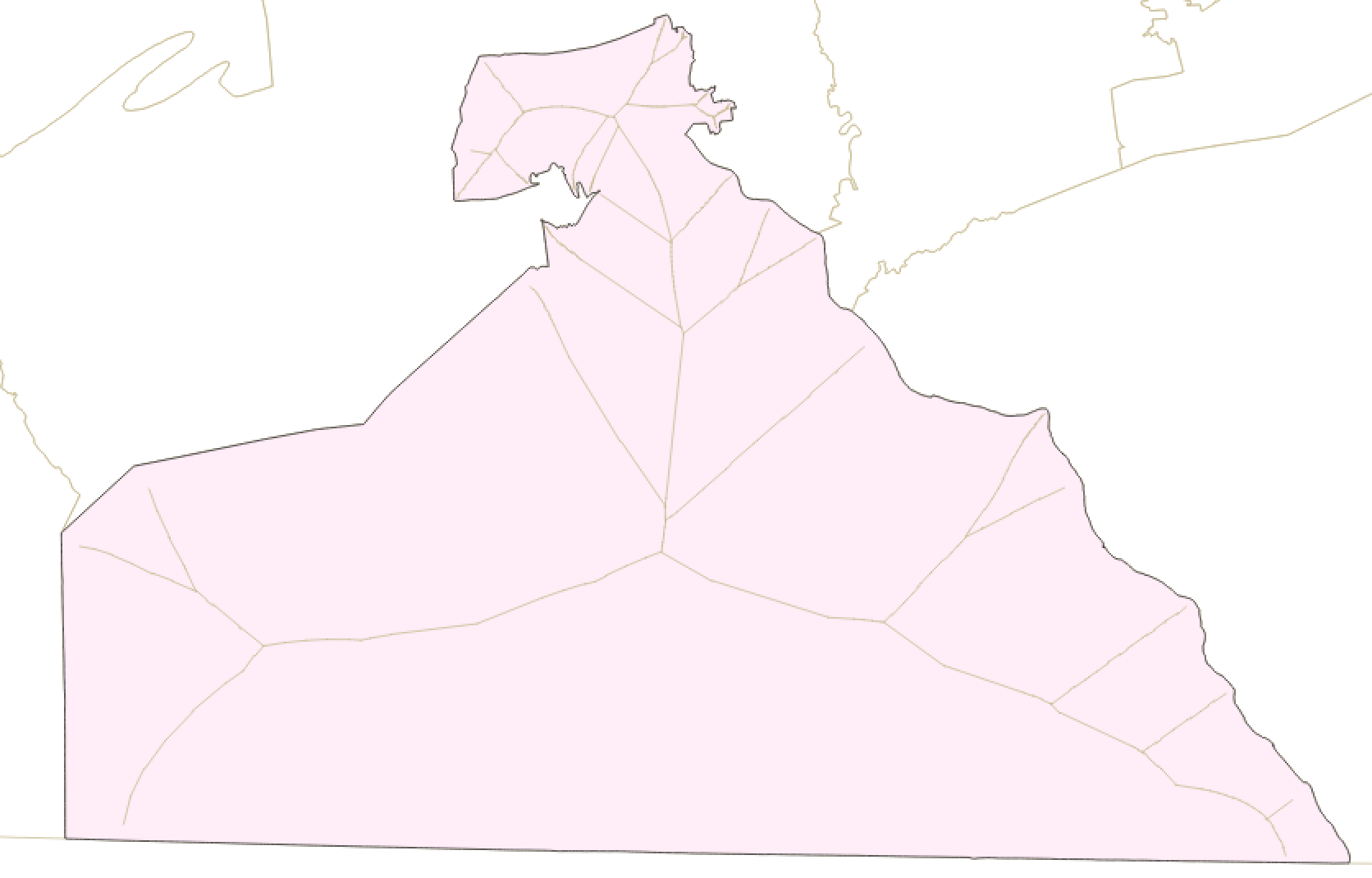}
\caption{Medial axes of PA 3 (left) and PA 4 (right)}\label{fig:PA03-04}
\end{figure}

In this district, which straddles the northwestern state corner, is somewhat boxy, with a few large, but somewhat irregular excisions in more rural territory. This yields a medial-hull ratio of 2.23, which places the district in Category 2 of gerrymandering evidence.

\subsubsection{Pennsylvania district 4}
Now, the largest district yet covers the eastern end of Pennsylvania’s south border, cresting up in triangular shape with two noticeable excisions along the way. Here we witness a possible issue as in Section \ref{subsec:noise} where noise on the edge of this large district causes lengthy medial limbs to stretch out from the more central spine. In doing so, the total medial length is bolstered, perhaps too much for what such noise should warrant in terms of gerrymandering. This noise is smoothed over for the hull axis, so the eventual medial-hull ratio of 2.14 is in some sense inflated, bringing the district into Category 2 of having been gerrymandered. While this claim does have merit in the northern and more populous tip, the district seems overall behaved in terms of shape.

\subsubsection{Pennsylvania district 5}

\begin{figure}
\includegraphics[width=0.475\textwidth]{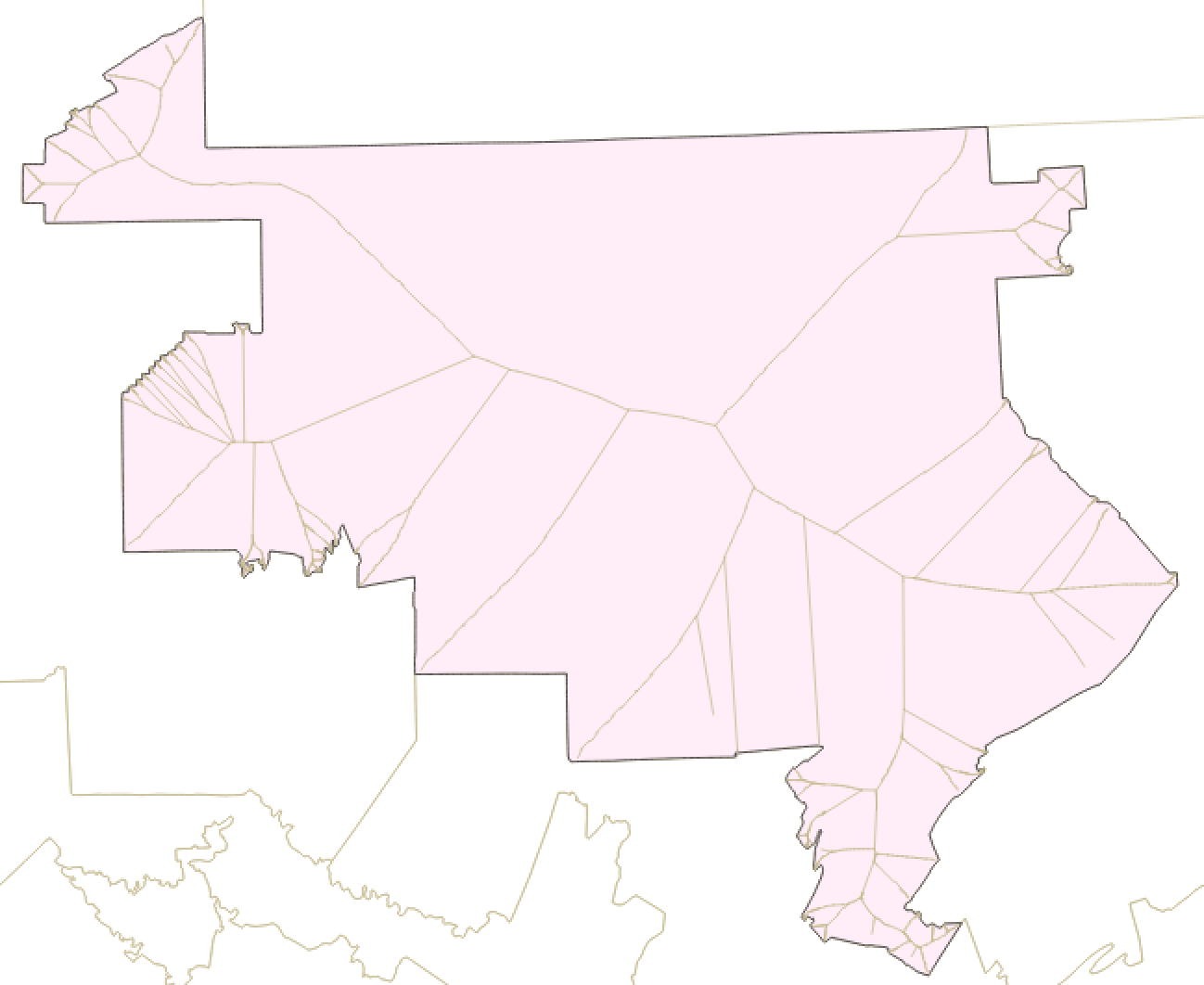}
\hspace{\fill}
\includegraphics[width=0.475\textwidth]{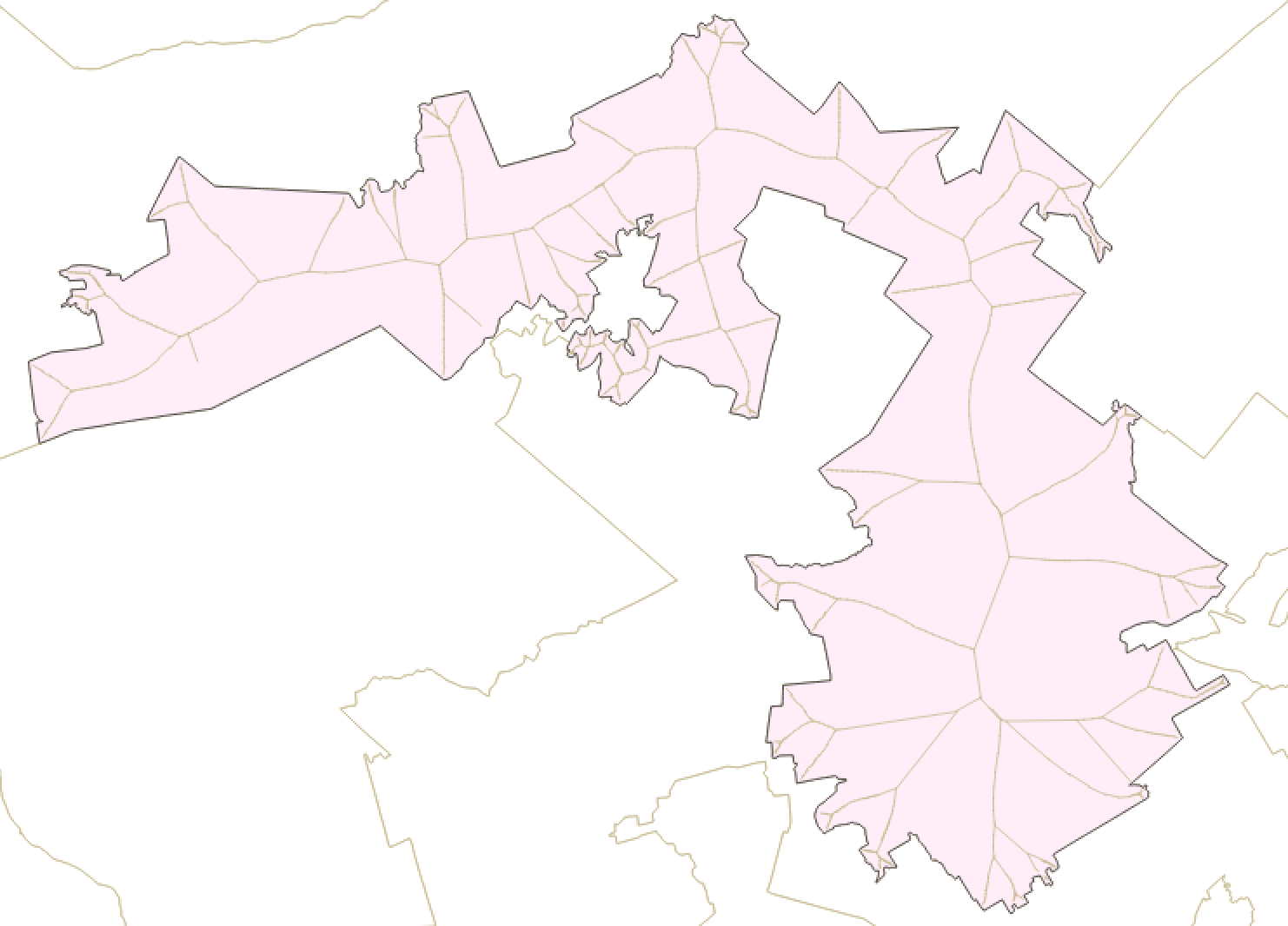}
\caption{Medial axes of PA 5 (left) and PA 6 (right)}\label{fig:PA05-06}
\end{figure}

The state’s largest district also has one of the highest medial-hull ratios: 3.10 (cat.~4). Certainly subject to bias against large area, this district witnesses long limbs in its medial axis that traverse a great distance from the central part of the district (and also of the medial axis) to the boundary where there is noise. Still, the southeastern tip and the western edge of the district seem too carved for comfort, and these certainly warrant further investigation into political intent.

\subsubsection{Pennsylvania district 6}
This district winds from just west of Philadelphia, to the north, then a seemingly unmotivated turn west. Several juts along the way are noticeable. Indeed, the winding shape of this district lands it nearly in Category 3 with a medial-hull ratio of 2.34 (cat.~2).

\subsubsection{Pennsylvania district 7}

\begin{figure}
\includegraphics[width=0.475\textwidth]{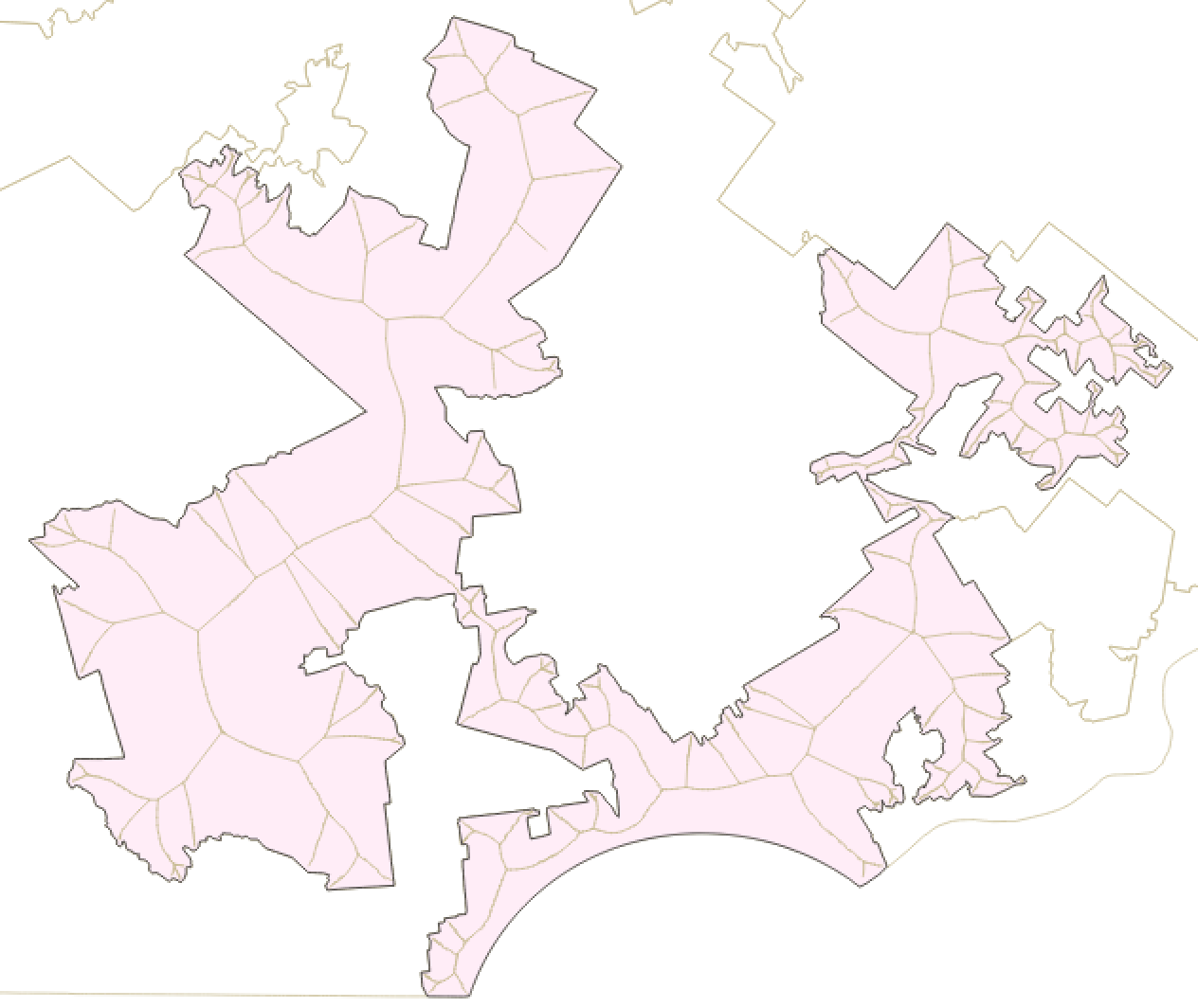}
\hspace{\fill}
\includegraphics[width=0.475\textwidth]{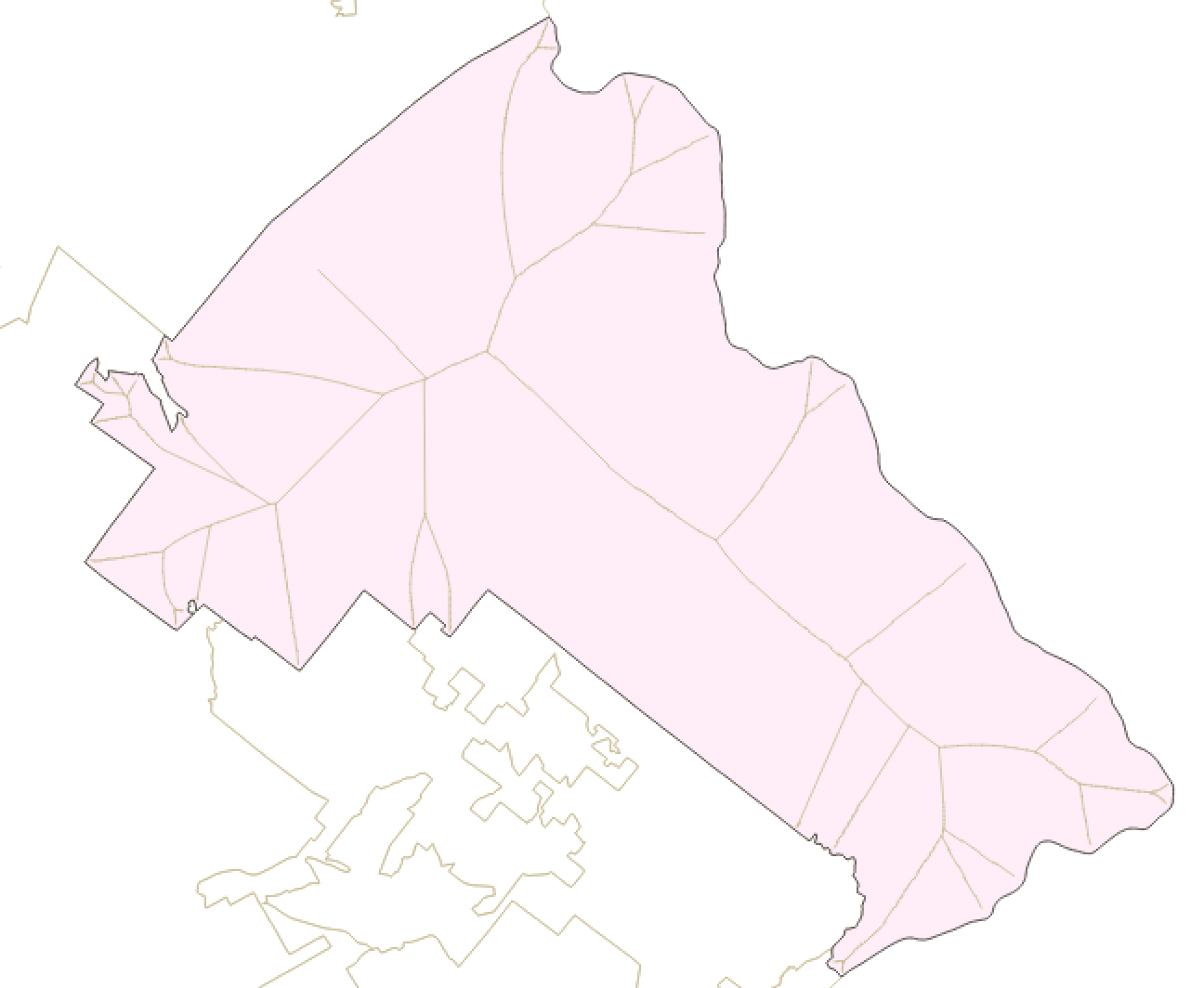}
\caption{Medial axes of PA 7 (left) and PA 8 (right)}\label{fig:PA07-08}
\end{figure}

Serving as the complement to Pennsylvania district 6 for much of its border, this district seems at first glance even more oddly shaped. Some have even gone so far as to nickname it ''Goofy Kicking Donald Duck." Indeed, the covered region stretches from north and west of Philadelphia to hugging the border Pennsylvania shares with Delaware, then winds around with major juts along the way. The medial-hull ratio of 2.47 lands the district in our Category 3 of gerrymandering.

\subsubsection{Pennsylvania district 8}
Tucked neatly against part of the eastern state border and otherwise fairly convex, this district bears the lowest medial-hull ratio of its state: 1.38 (cat.~1). We witness an Italy-shaped jut out of the western side of the district, but otherwise see no alarm in the district's shape.

\subsubsection{Pennsylvania district 9}

\begin{figure}
\includegraphics[width=0.475\textwidth]{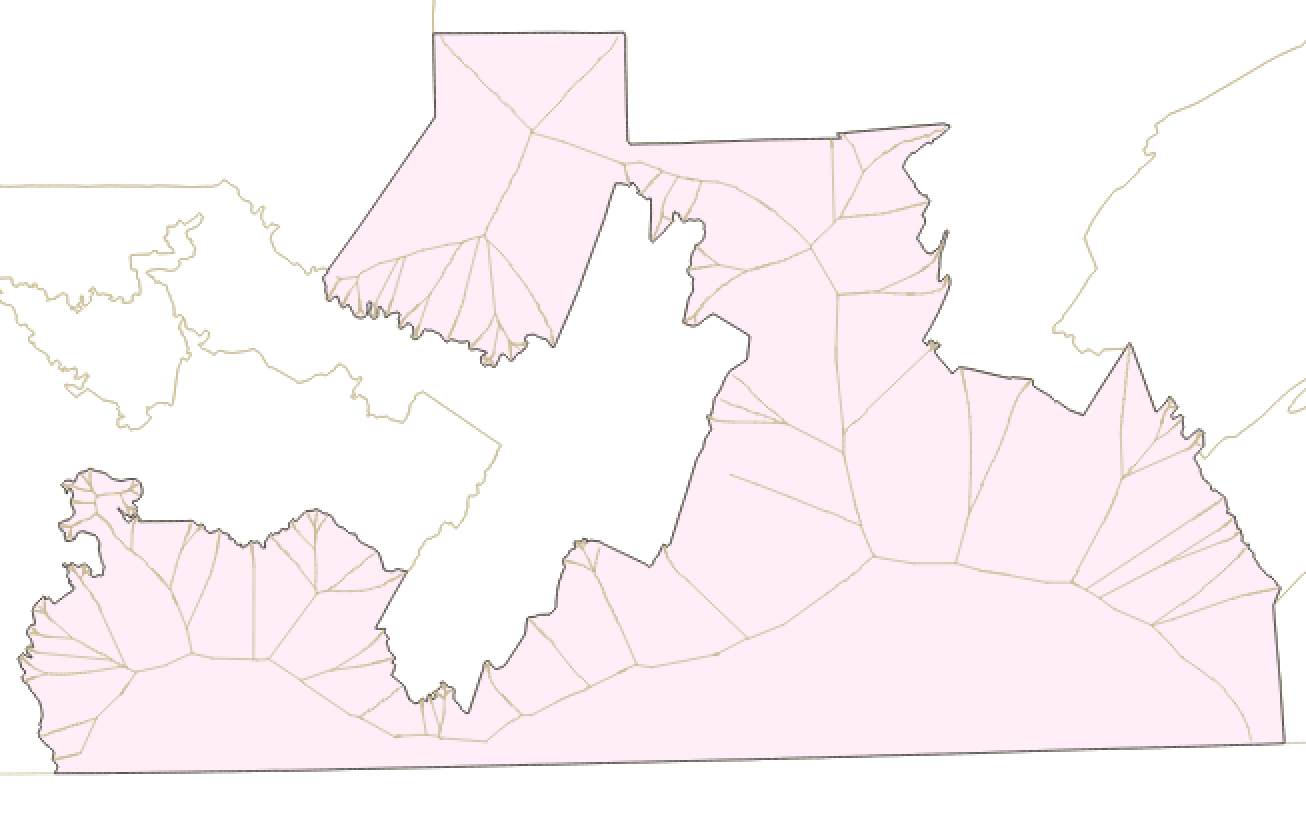}
\hspace{\fill}
\includegraphics[width=0.475\textwidth]{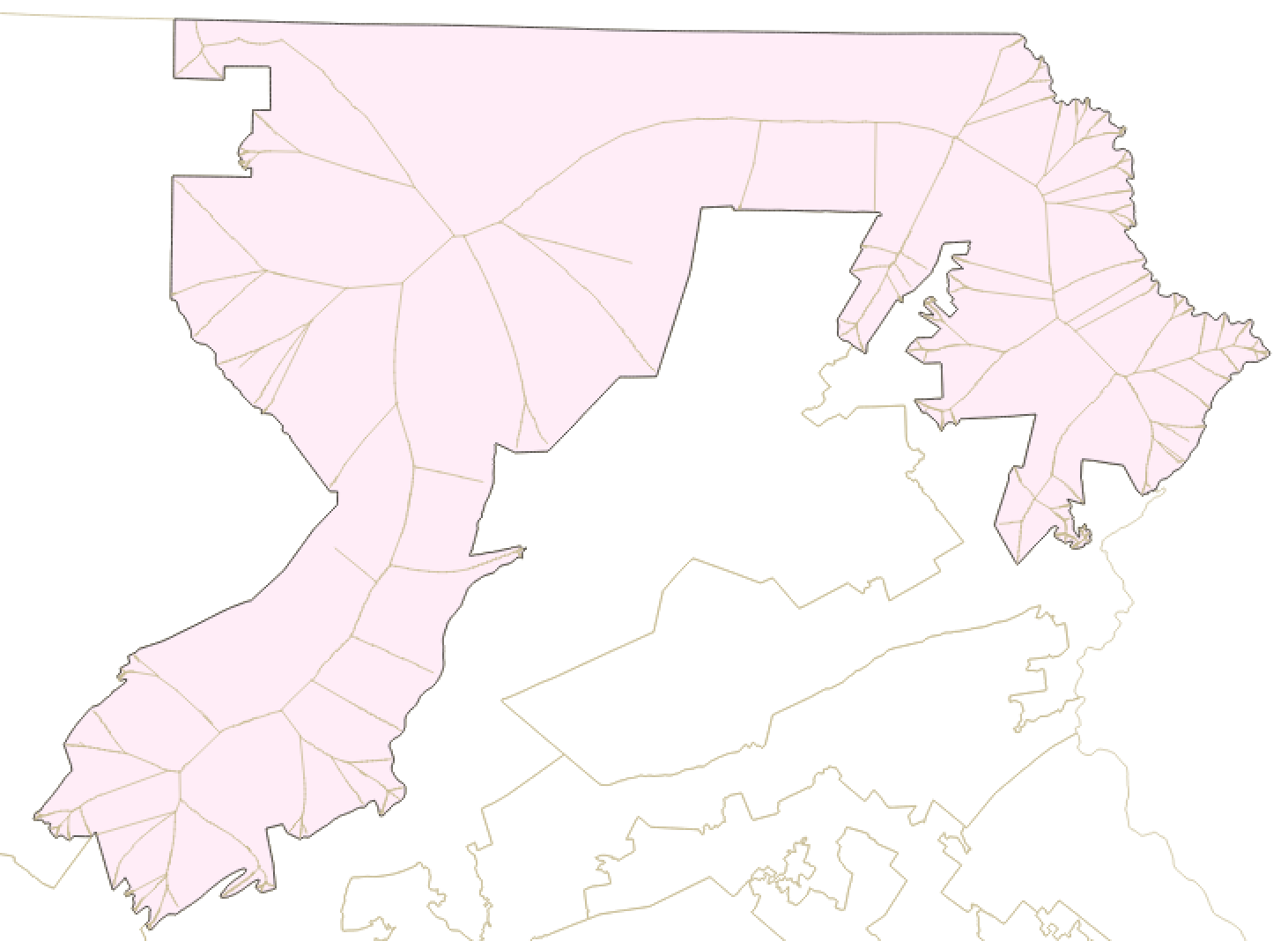}
\caption{Medial axes of PA 9 (left) and PA 10 (right)}\label{fig:PA09-10}
\end{figure}

This district is the third-largest by area and bears the second-highest medial-hull ratio with 3.35 (cat.~4). It sits along the western side of Pennsylvania's southern border and stretches toward the north. However, a massive excision from the west occurs on the border it shares with PA 12. This region is the source of much of the district's medial axis, though it is smoothed over with the convex hull. Thus, the medial-hull ratio is high, and we are near-certain there is gerrymandering at play.

\subsubsection{Pennsylvania district 10}
Hugging the northeast corner of the state, this district falls into low Category 2 of our gerrymandering gauge with a medial-hull ratio of 2.01. While not ostensibly offensive, the shape does bear a ``U" pattern and have some noticeable juts complementing PA 17. We also remark that this district is the second largest in the entire state, so it may suffer from bias against larger districts via our model parameters.

\subsubsection{Pennsylvania district 11}

\begin{figure}
\includegraphics[width=0.475\textwidth]{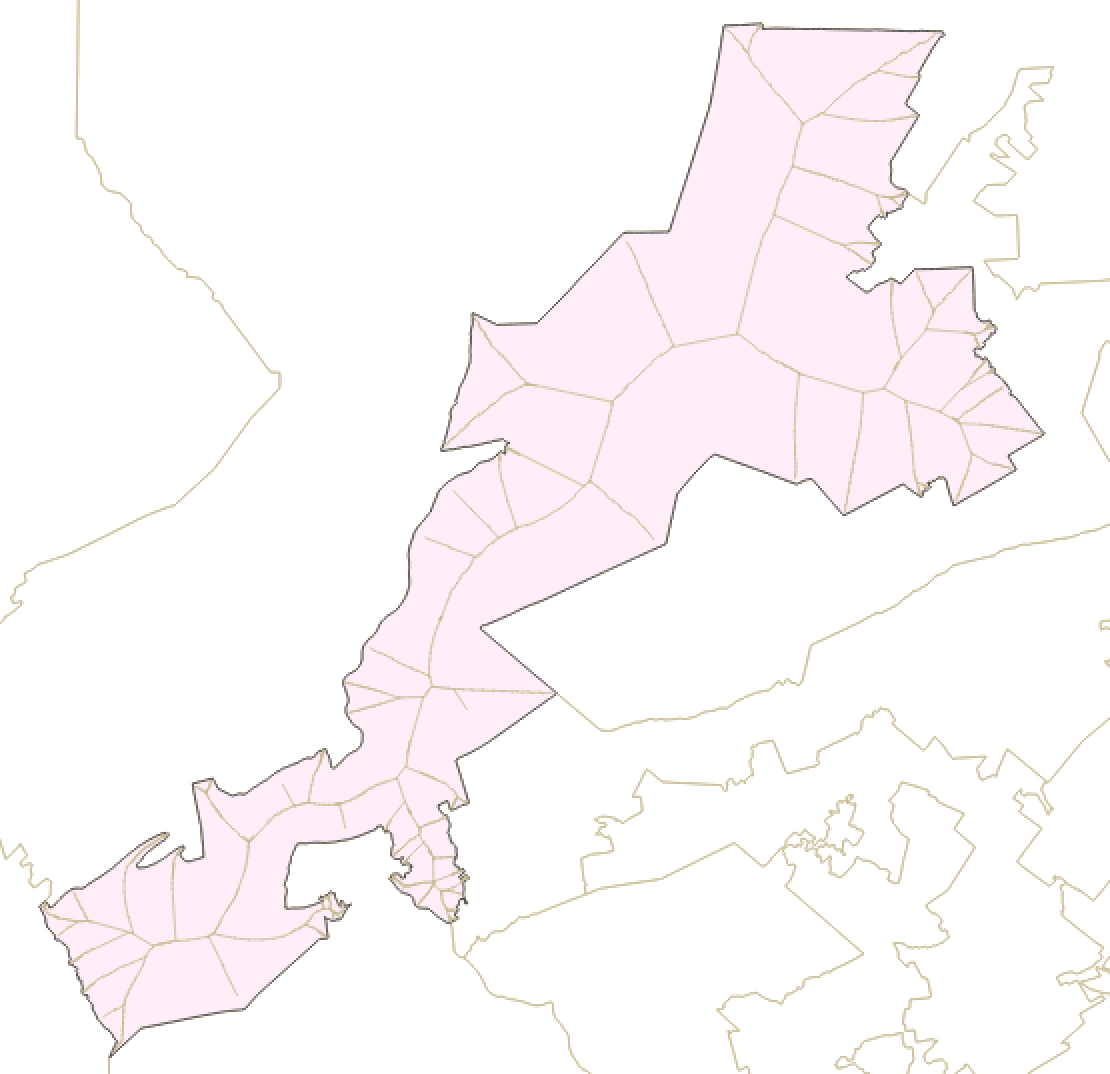}
\hspace{\fill}
\includegraphics[width=0.475\textwidth]{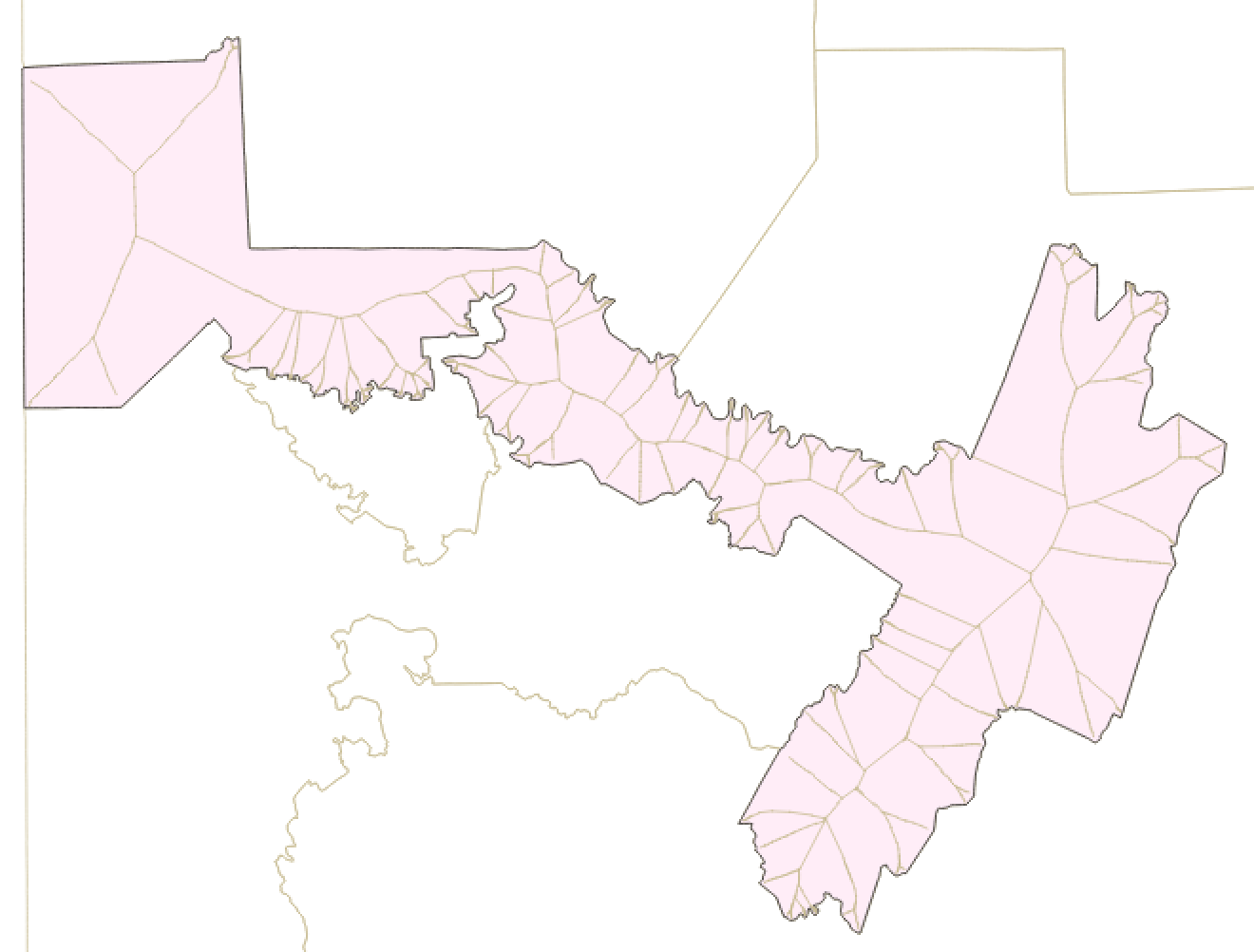}
\caption{Medial axes of PA 11 (left) and PA 12 (right)}\label{fig:PA11-12}
\end{figure}

District 11 traverses swaths of land near Harrisburg up to Scranton. The high medial-hull ratio of 3.20 (cat.~4) indicates strong evidence of gerrymandering. However, the district's medial length is largely derived from noise (which the convex hull smooths over) along its somewhat long, rectangular shape. This bolsters the ratio perhaps a bit too much, but overall demonstrates a successful instance of penalizing districts which are much longer than they are thick.

\subsubsection{Pennsylvania district 12}
Here lies much of the demise from PA 9. The western portion of this district seems innocent enough, though it is a long and skinny rectangle. The eastern half, however, punches a sizable chunk out of PA 9, heightening each of these districts' medial-hull ratios. This district wields a ratio of 3.25, pointing to considerable evidence of gerrymandering in Category 4.

\subsubsection{Pennsylvania district 13}

\begin{figure}
\includegraphics[width=0.475\textwidth]{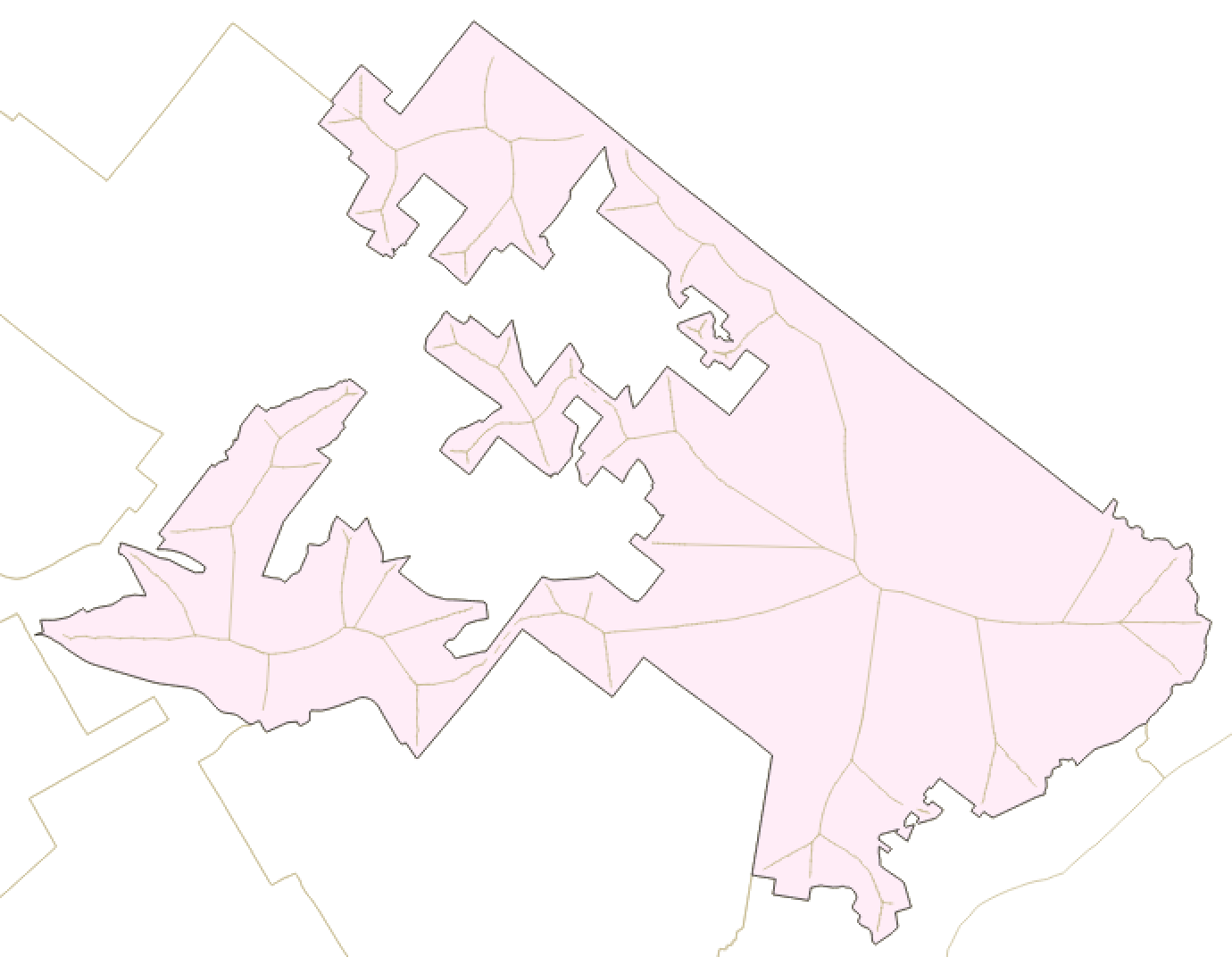}
\hspace{\fill}
\includegraphics[width=0.475\textwidth]{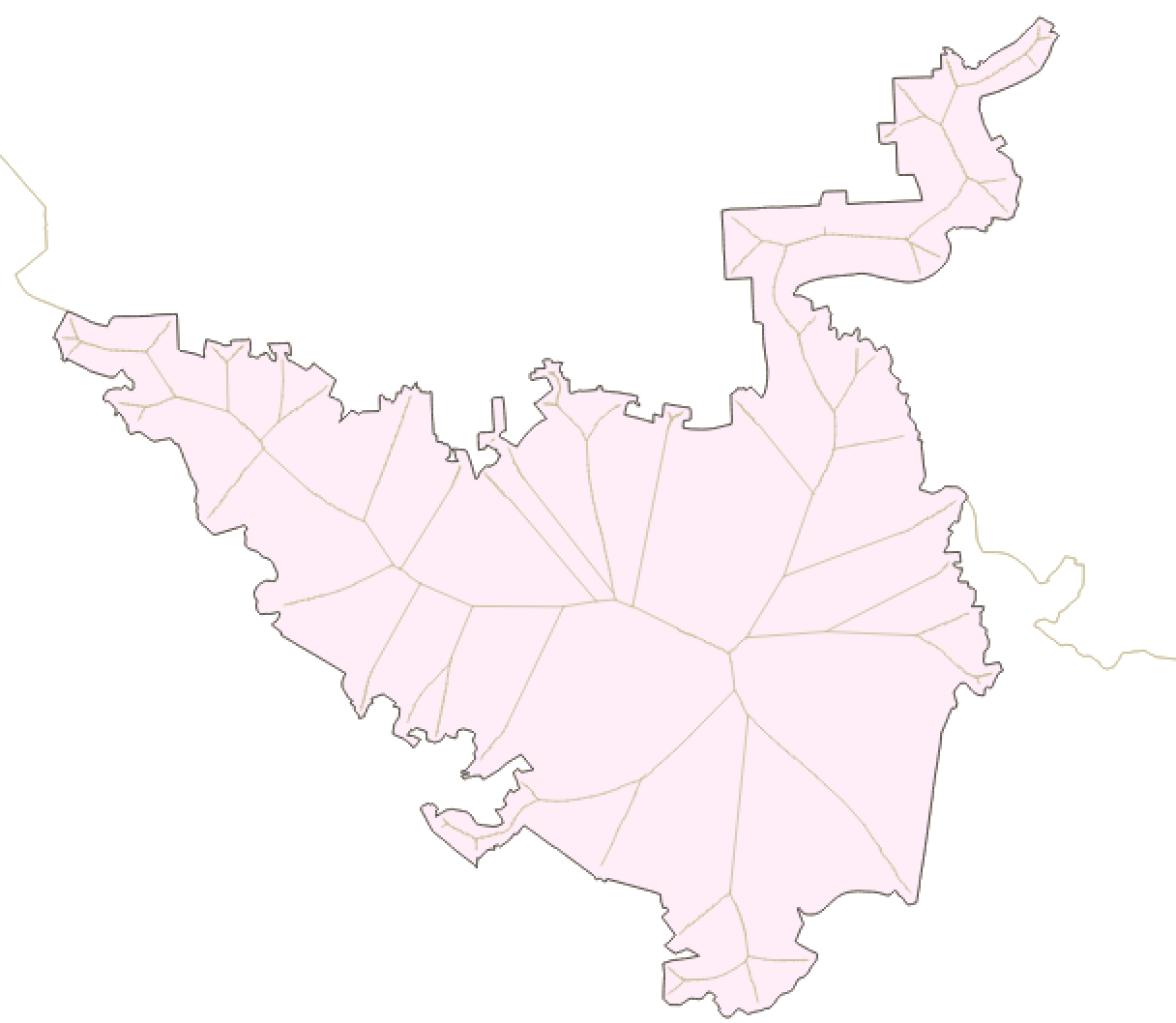}
\caption{Medial axes of PA 13 (left) and PA 14 (right)}\label{fig:PA13-14}
\end{figure}

This district lies near the eastern state border, covering much of the land north-northeast of Philadelphia but possessing clear tracts of land excised from it and passed off to PA 7. Agreeably, its medial-hull ratio of 2.80 places it just barely in Category 4 of our gerrymandering accusation.

\subsubsection{Pennsylvania district 14}
Largely covering Pittsburgh, this district resembles a goose in shape. With a medial-hull ratio of 2.65, the district falls within Category 3 of gerrymandering indication. Indeed, while the district's shape is somewhat ovular, it has an unusual winding extension out to the northeast, which hints at malicious political intent.

\subsubsection{Pennsylvania district 15}

\begin{figure}
\includegraphics[width=0.475\textwidth]{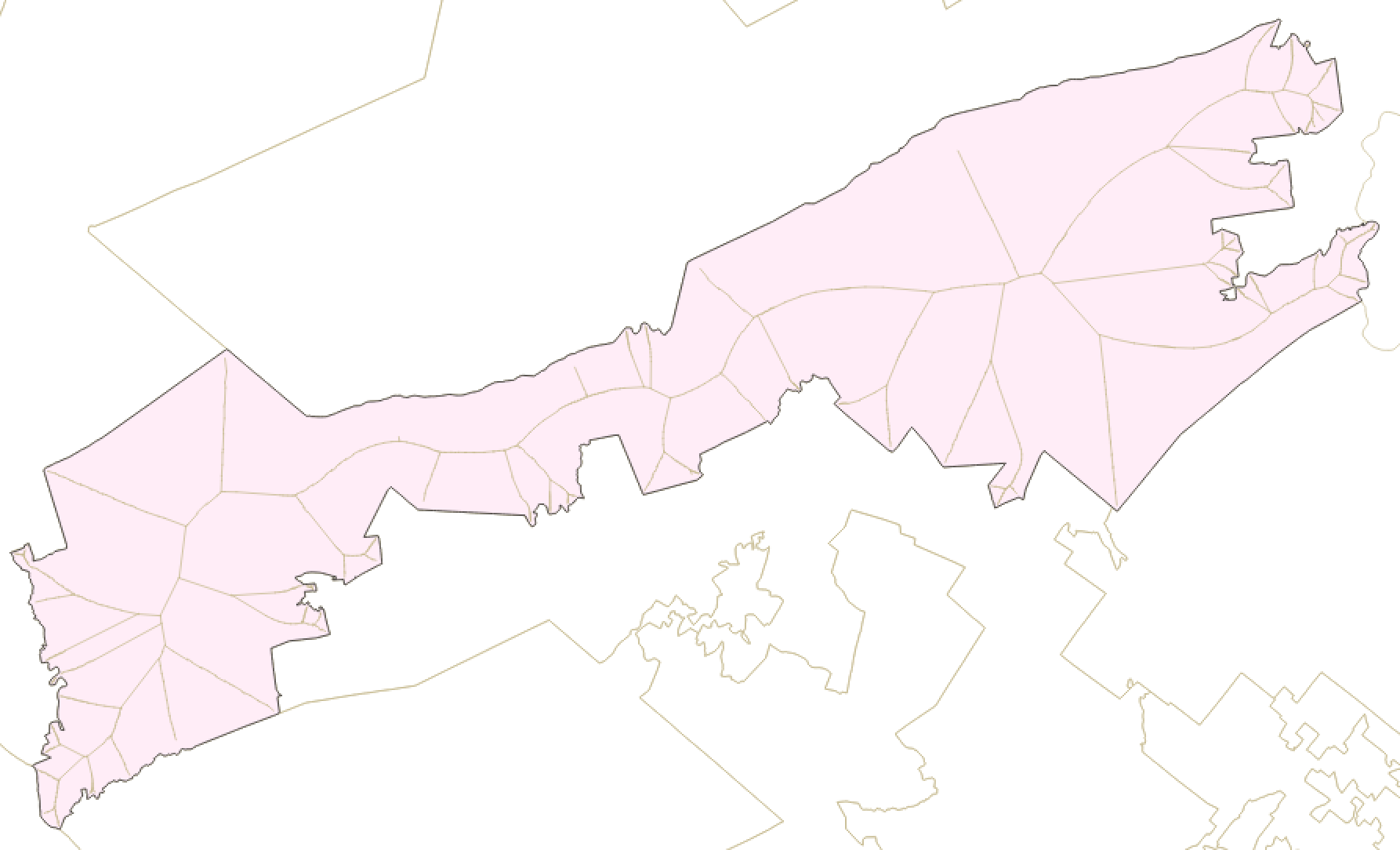}
\hspace{\fill}
\includegraphics[width=0.475\textwidth]{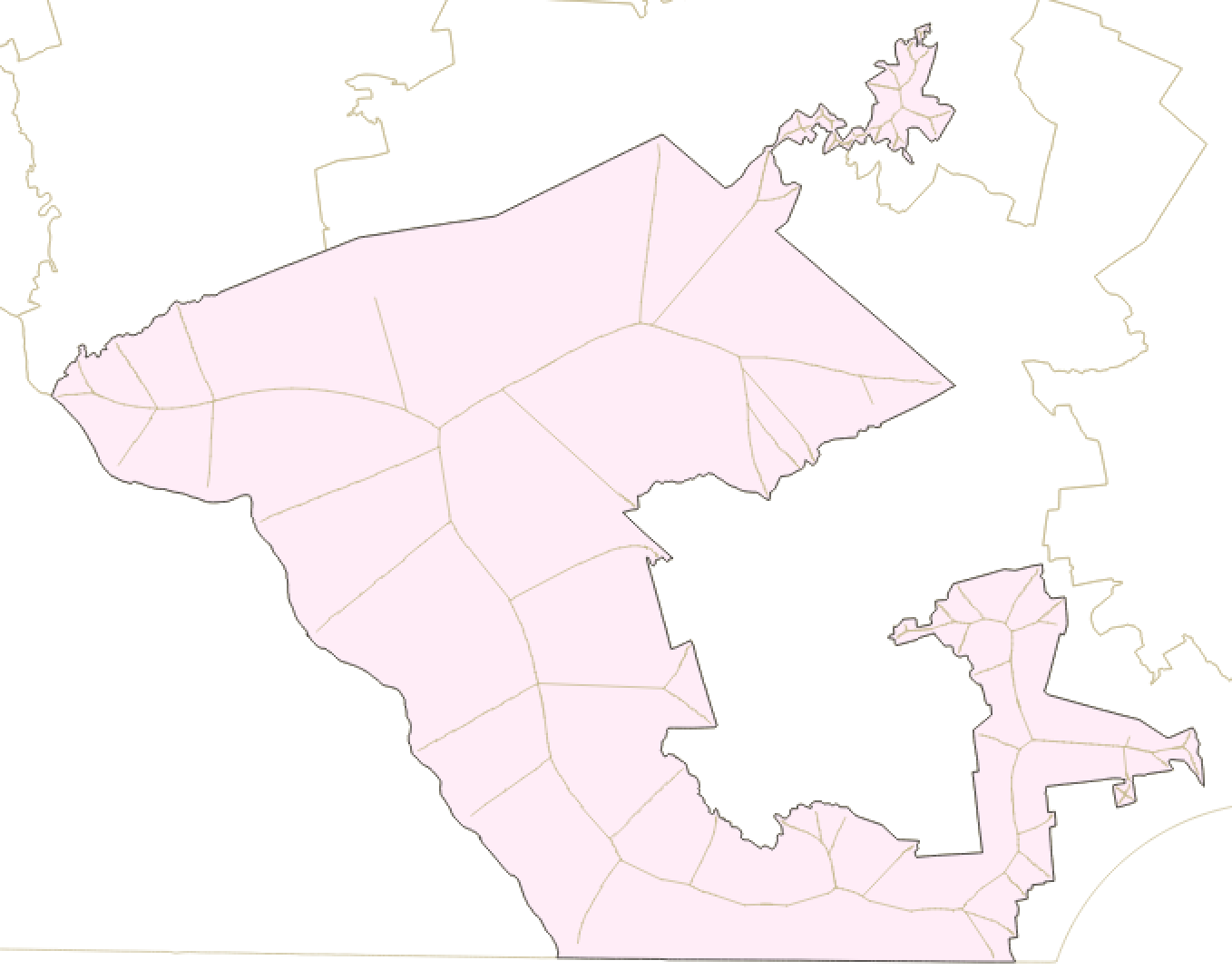}
\caption{Medial axes of PA 15 (left) and PA 16 (right)}\label{fig:PA15-16}
\end{figure}

The shape of this district resembles that of PA 11, albeit with fewer and smaller juts. This time, the medial-hull ratio is 2.48 (cat.~3), which seems reasonable given its long, thin stretch of territory and the inward juts on its eastern end.

\subsubsection{Pennsylvania district 16}
On one side forming the complement to PA 7, this Mexico-shaped district has a single, intense jut out, but otherwise winds modestly. Its medial-hull ratio of 2.71 (cat.~3) confirms that it bears probable evidence of gerrymandering.

\subsubsection{Pennsylvania district 17}

\begin{figure}
\includegraphics[width=0.475\textwidth]{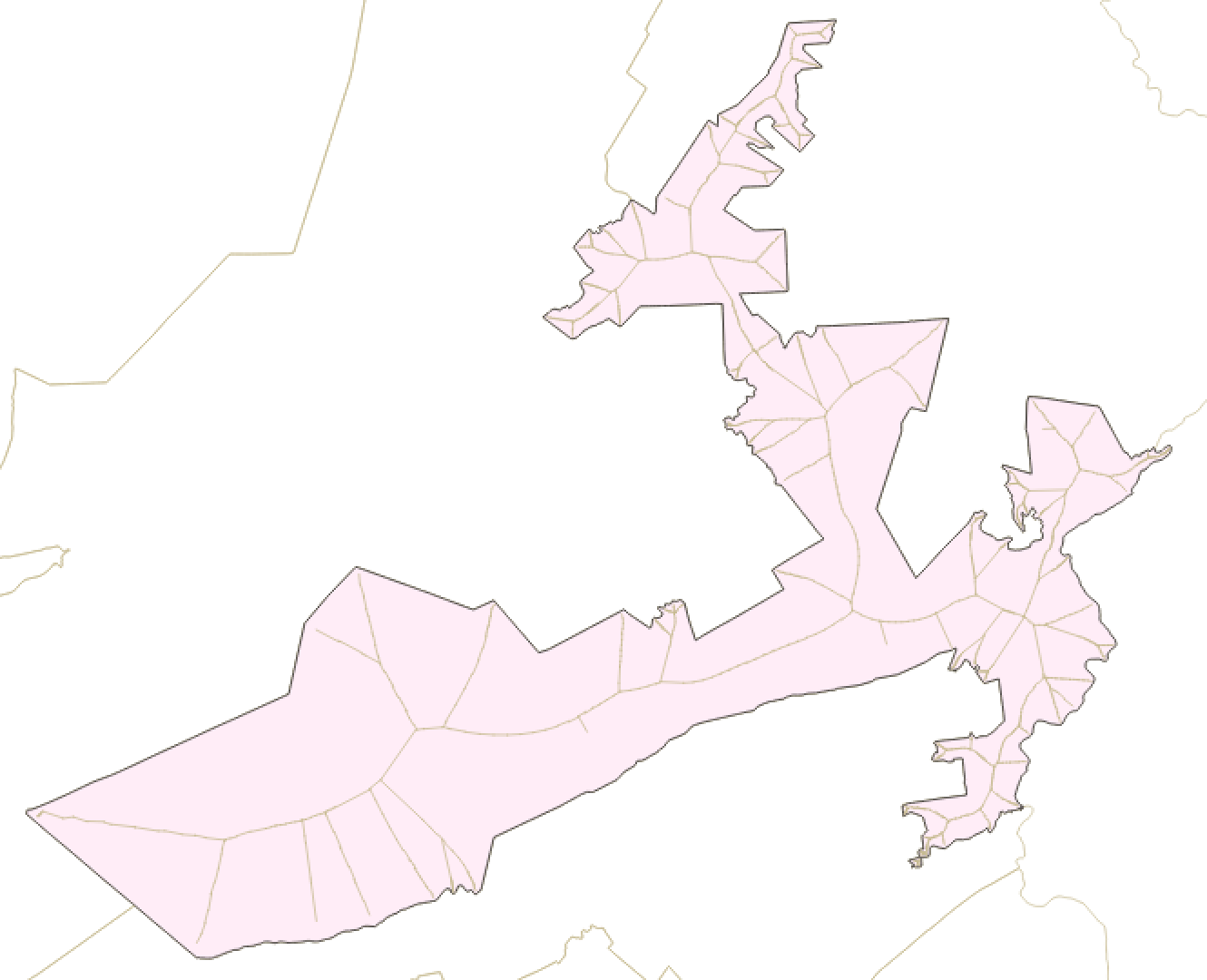}
\hspace{\fill}
\includegraphics[width=0.475\textwidth]{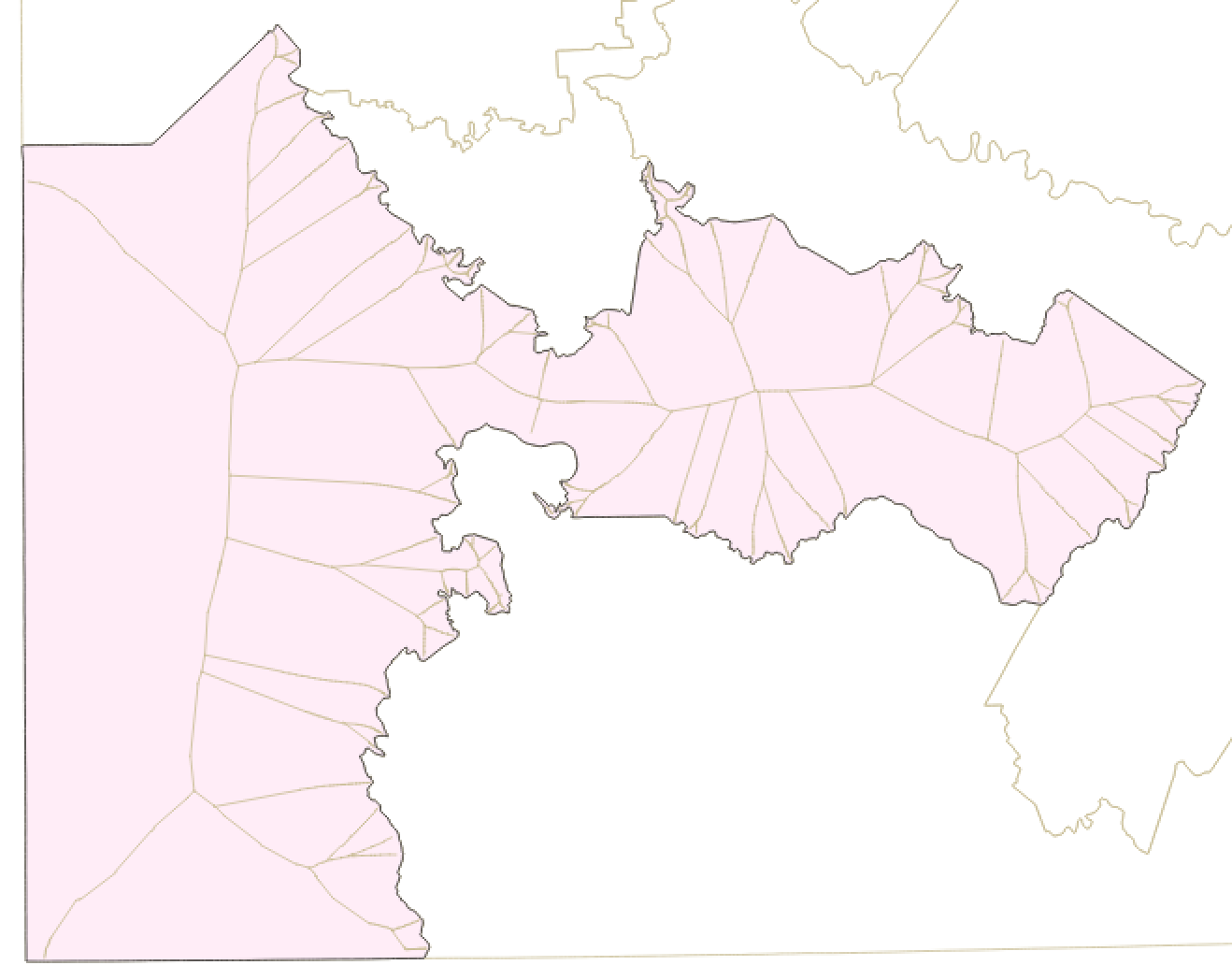}
\caption{Medial axes of PA 17 (left) and PA 18 (right)}\label{fig:PA17-18}
\end{figure}

Half of this district seems boxy and pleasingly rectangular, while the other half is more errant, splicing into two noisy tendrils. Unfortunately, this is an example of a district which has a well-behaved portion masking the evidence of gerrymandering from a wilder other half under our ratio metric. The medial-hull ratio of this district almost, but not quite, qualifies it for Category 2. The district's shape behavior on the west side provides a buffer for the carving done on the east end, bringing down its overall ratio and steadying it close to 2.00. We remark that had the medial-hull ratio been calculated on just the eastern portion of the district, subsections of the district gave way to ratios as high as 2.50, flagging probable evidence of gerrymandering with Category 3.

\subsubsection{Pennsylvania district 18}
This district bears the largest medial-hull ratio of any in the state, with 3.50. Sitting in the southwestern corner of Pennsylvania, the district nearly pinches off completely, then juts out a large plot of land from its eastern edge. This is almost certainly evidence of gerrymandering, as corroborated by its qualification to Category 4 of our gerrymandering analysis.

\subsection{Pennsylvania (remedial)}

\begin{figure}
\centerline{\includegraphics[height=2.5in]{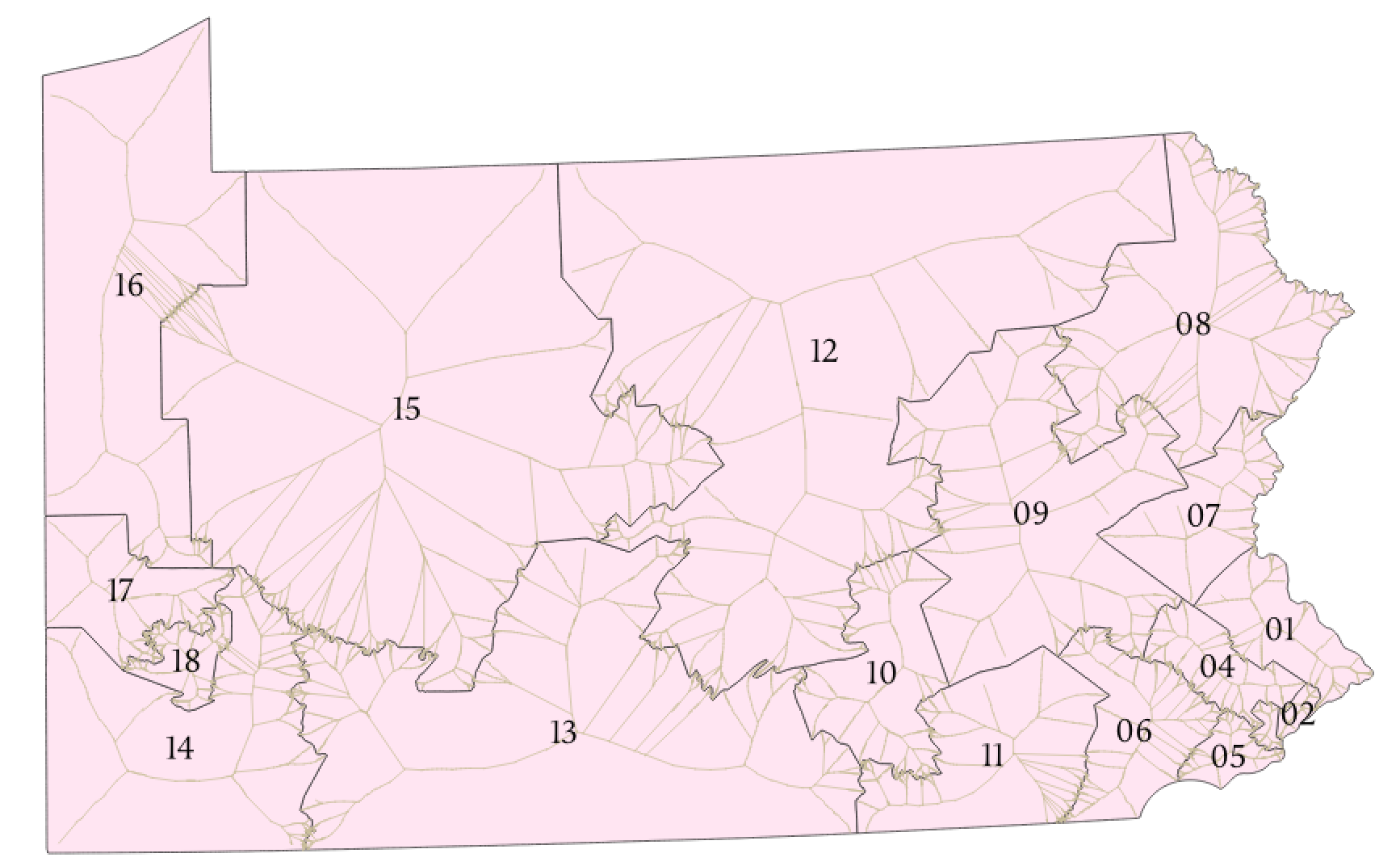}}
\caption{Medial axes of remedial Pennsylvania (2018) districts}\label{fig:PAremmed}
\end{figure}

We now discuss the Pennsylvania remedial map and how it improves the state's performance under our gerrymandering analysis.

\begin{figure}
\centerline{\includegraphics[height=2.5in]{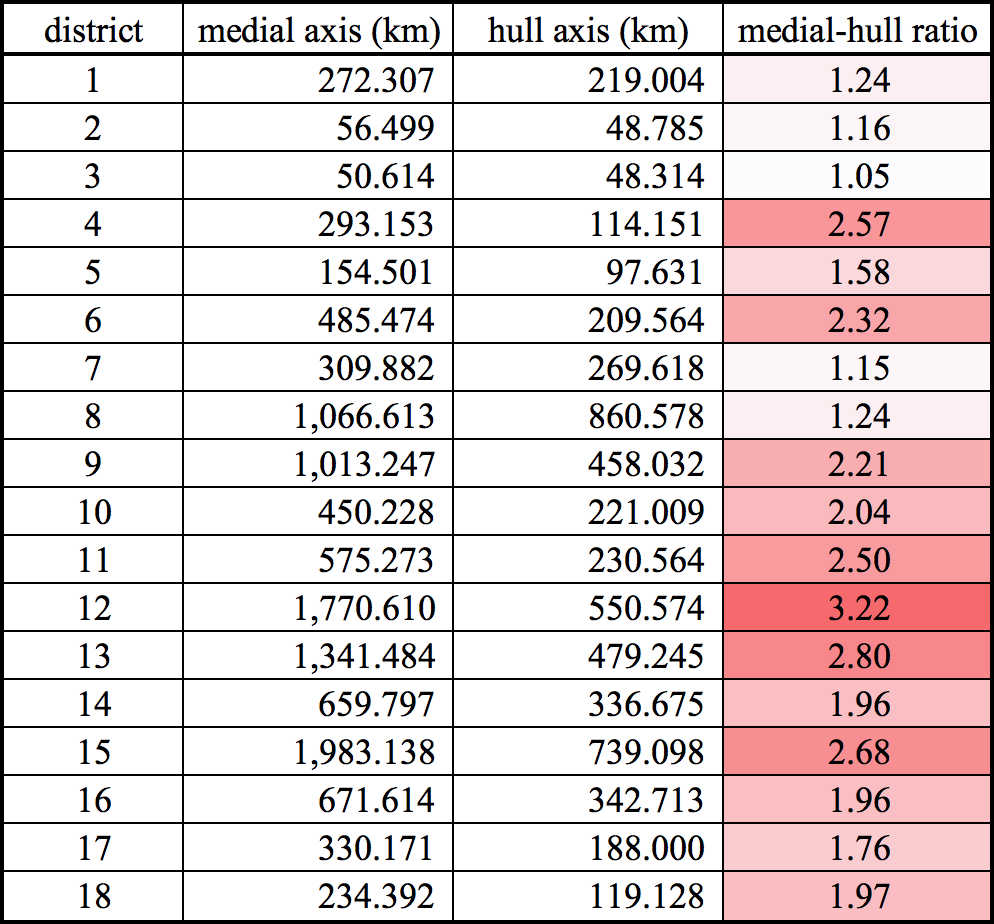}}
\caption{\label{fig:PAnew} Pennsylvania remedial district data}
\end{figure}

\begin{remark}
Since the remedial map represents an extreme reworking of the district shapes and arrangement, district numbers are not meant to correspond. For example, the old district 1 should not be directly compared to the remedial district 1, as they are nearly disjoint and not intended to align. Instead, we will compare the maps in whole and trends in district performance under our medial-hull ratio.
\end{remark}

When viewing the map, we immediately see more ostensibly compact shapes. In the eastern portion of the state, we witness larger juts, if any, and generally rectangular shapes. Pennsylvania remedial districts 1, 2, 3, 7, and 8 all achieve ratios less than 1.25. Generally, they have at most one worrisome jut (aside from the state border) and overall seem fair. Remedial districts 4 and 6 land in low Category 3 and high Category 2 of gerrymandering evidence, respectively. Indeed, much of their superfluous medial length stems from the noise on their shared border. Other than this shared edge and another large jut each, however, these districts behave well.

In the north-central region of the state, remedial district 12 bears the highest medial-hull ratio with 3.22. Indeed, the region is large (so it may suffer partially from bias against large districts), but also has an errant and noisy southeast side. This certainly seems to imply political motivation, as confirmed by the Category 4 gerrymandering qualification of the ratio. Similarly, remedial districts 13 and 15 suffer from the two next-highest medial-hull ratios of 2.80 (cat.~4) and 2.68 (cat.~3), respectively. Together, these complement remedial PA 12 and thus share those errant boundaries on their east sides. In the northwest corner of district 13, we also observe an unusually noisy jut toward Pittsburgh. The other remedial districts in the western half of the state raise little alarm, whether by visual inspection or our medial-hull ratio.

\subsubsection{Improvements over original map}
Now, we compare the old Congressional district map to its remedial version with celebration. The statewide average fell from 2.50 (cat.~3) to 1.97 (cat.~1). The number of Category 4 districts dropped from 6 to 2, Categories 3 and 2 districts each dropped from 4 to 3 in count, and the coterie of alarm-free districts in Category 1 happily expanded from 4 to 10! Indeed, according to our medial-hull ratio, the remedial mapping largely rid the state of gerrymandering issues.

\subsection{North Carolina}

\begin{figure}
\centerline{\includegraphics[height=2.5in]{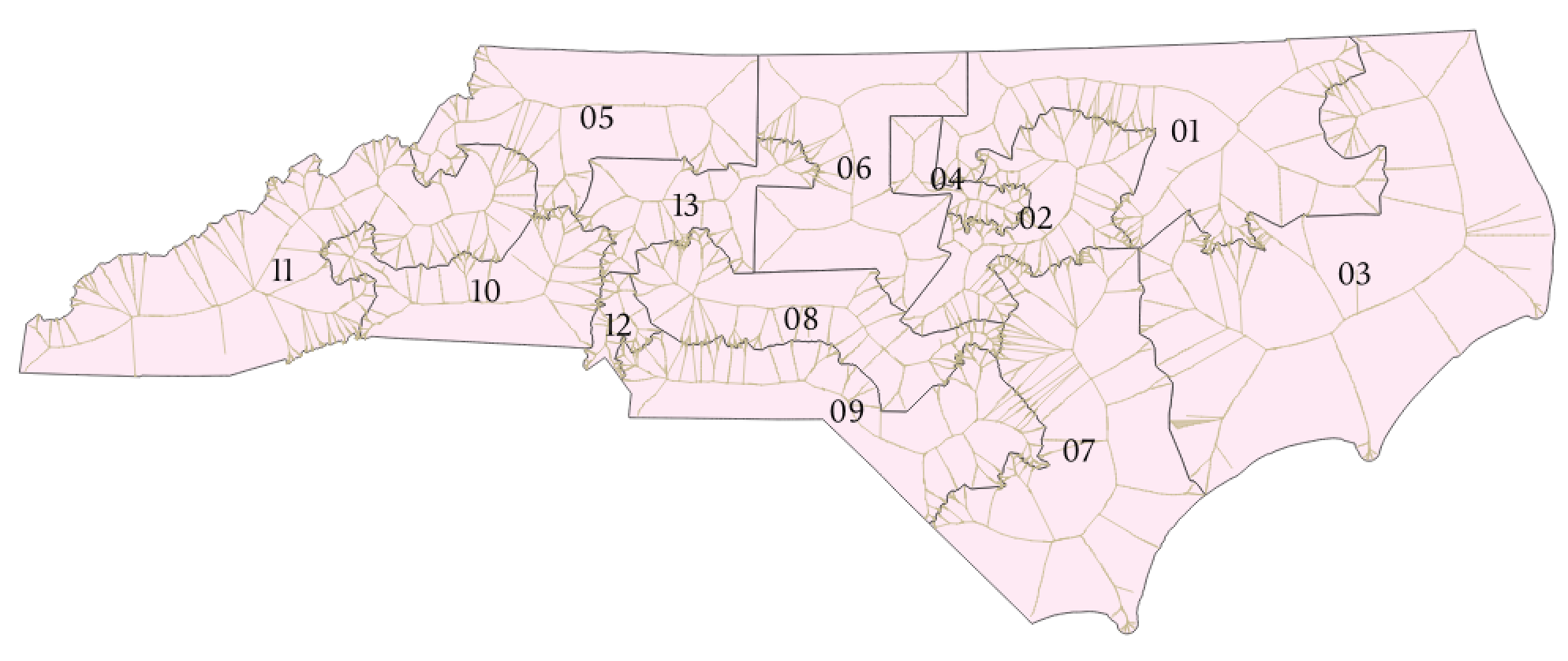}}
\caption{Medial axes of North Carolina districts}\label{fig:NCmed}
\end{figure}

We take a closer look at remarkable districts in North Carolina\footnote{State FIPS code 37}, grouped by their category rankings under the medial-hull ratio. Overall, this state possesses 2 districts in Category 4, 5 in Category 3, 1 in Category 2, and 5 in Category 1. These compose a statewide average of 2.32 (cat.~2).

\begin{figure}
\centerline{\includegraphics[height=2.0in]{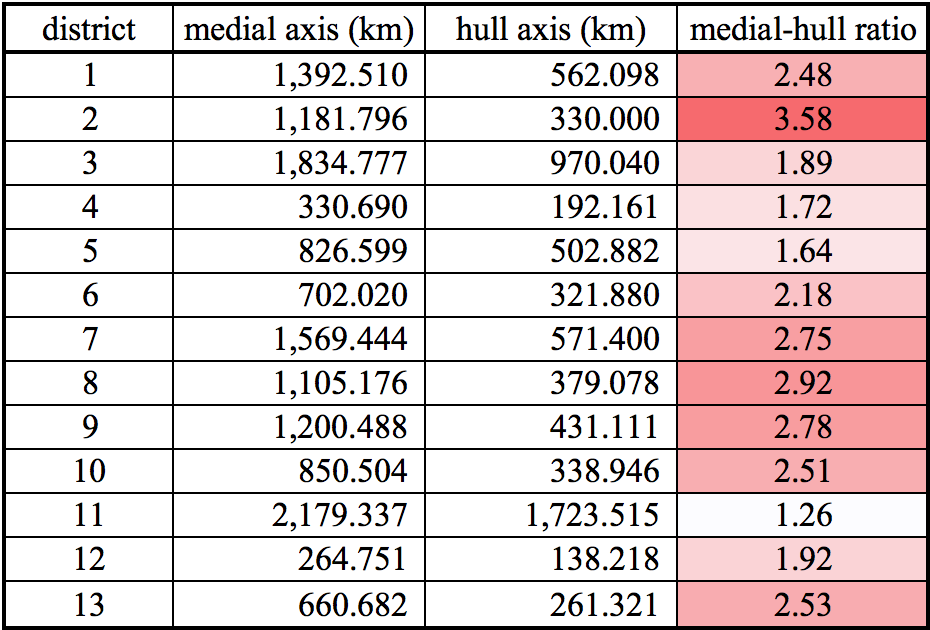}}
\caption{\label{fig:NC} North Carolina district data}
\end{figure}

\subsubsection{Category 4}
District 2 bears a statewide record medial-hull ratio of 3.58. Indeed, its boundary bears considerable noise throughout. Further, it follows the shape of a backwards `S', winding around itself and other sizable juts from other districts. Together, these characteristics suggest that gerrymandering is at play.

District 8 has the state's second-highest medial-hull ratio with 2.92. Its shape misbehaves in a different way from NC 2. Instead of winding around, it stretches horizontally and curls at either end. Moderate noise along most of its border also contributes heavily to its total medial length. This district spanning with such breadth but little width raises flags to malicious political intent.

\subsubsection{Category 3}
District 1 falls cleanly within Category 3 with a medial-hull ratio of 2.48. Hugging the state's northern border, its western end suddenly plunges southward and its eastern end juts down before splitting into three large regions. While these juts are large, they seem to mostly be accommodating the errant presence of NC 2, as discussed above. Thus, this district is on high watch for evidence of gerrymandering, but not itself the most offensive.

District 9 nearly qualifies for Category 4 with a medial-hull ratio of 2.78. Indeed, this district appears to mimic the shape of NC 8, which sits directly atop it. Following the southern North Carolina border, district 9 spans a great horizontal distance while remaining relatively thin, then curving upward at its eastern tip. This likely falls victim to the same gerrymandering that may have motivated the drawing of NC 8.

\subsubsection{Category 2}
District 6 is the only North Carolina district in this category. Indeed, it bears a shape similar to that of the capital letter `I' with a curved southeastern tail. While the jut into its western side also raises slight alarm, the district is otherwise quite boxy and compact, justifying its placement in Category 2 for gerrymandering evidence.

\subsubsection{Category 1}
District 4 seems to be the only North Carolina district in Category 1 that perhaps should be further inspected. Composed of two boxy regions, this district has a central, thin bridge that holds it together. This pattern does not seem intuitively compact, yet remains unpenalized by our medial process, which does not account for thickness in central regions of each district.

District 11 is nestled tightly into the state's southwestern corner. This mounts a good bit of medial length both onto the district and its clipped convex hull, thus neutering its medial-hull ratio. Indeed, the district has a noticeable jut inward from NC 10 (cat.~3), but otherwise behaves as well as it can, subject to that state border. Thus, our medial-hull ratio metric excuses this and does not flag it for gerrymandering evidence.

\subsection{Florida}

\begin{figure}
\centerline{\includegraphics[height=2.5in]{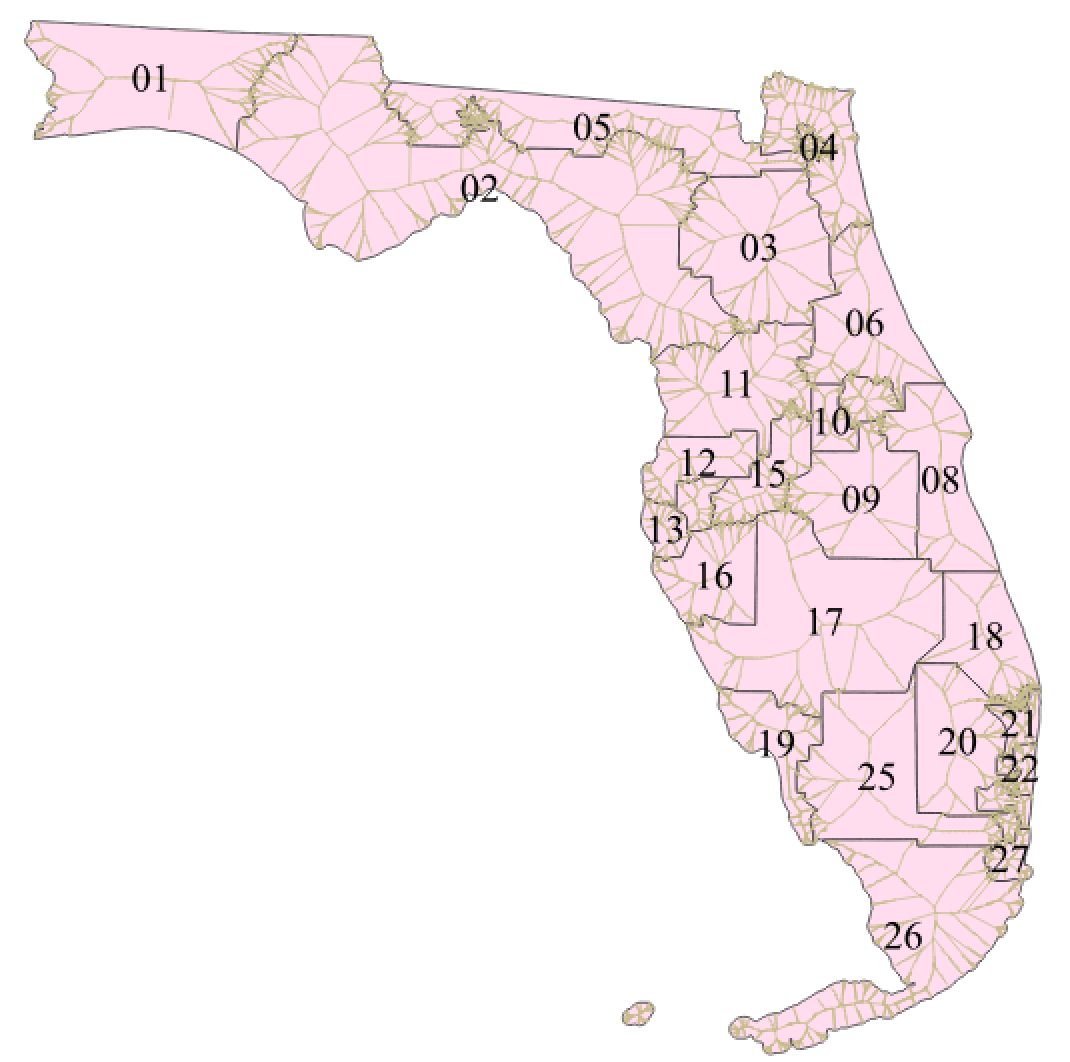}}
\caption{Medial axes of Florida districts}\label{fig:FLmed}
\end{figure}

We take a closer look at remarkable districts in Florida\footnote{State FIPS code 12}, grouped by their category rankings under the medial-hull ratio. Overall, this state possesses no districts in Category 4, 3 each in Category 3 and Category 2, and 21 in Category 1. These compose a statewide average of 1.65 (cat.~1).

\begin{figure}
\centerline{\includegraphics[height=3.5in]{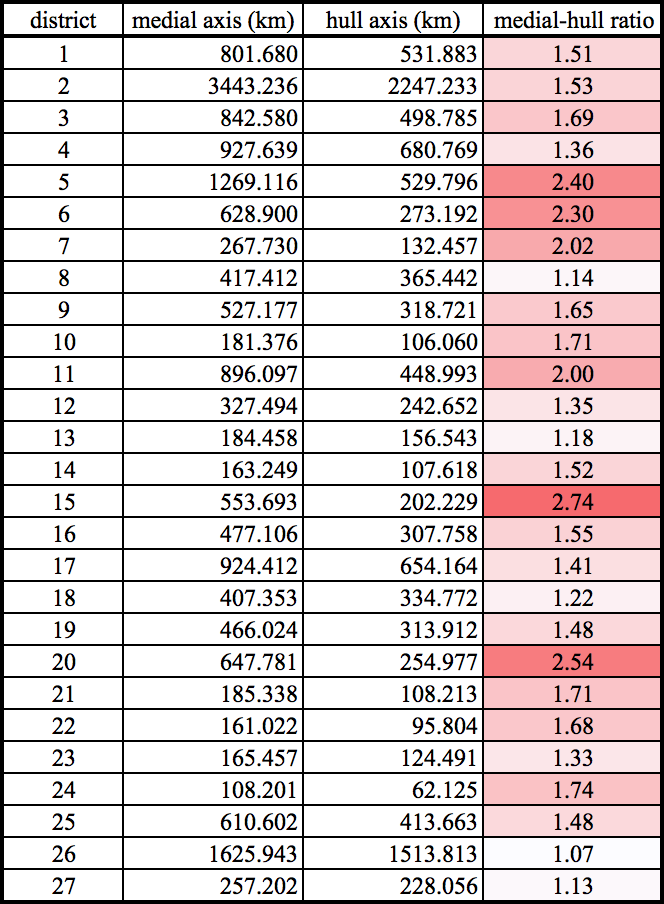}}
\caption{\label{fig:FL} Florida district data}
\end{figure}

\subsubsection{Category 3}
District 15, taking the crown for Florida's highest medial-hull ratio (2.74), has a thin, nearly ovular shape. Considerable noise along much of its border adds to the medial length and alerts us to evidence of gerrymandering.

The west side of district 20 behaves well. However, a long and spindly arm extends from its northeast side, both signaling gerrymandering and considerably increasing the district's medial-hull ratio to 2.54. Another arm due south helps solidify this district in Category 3.

The longest and thinnest rectangular district witnessed yet, district 5 straddles the north border of Florida. However, this is not enough to cover the effect of an inward jut on the west or a jagged arm extending from the east side of the district. These certainly justify the district's placement in this category. 

\subsubsection{Category 2}

District 6, with a medial-hull ratio of 2.30, shares a noisy border with each FL 7 and FL 11. It also sees a peculiar jut with FL 3. While remaining otherwise compact, between these three red flags, district 6 deserves its place solidly within Category 2.

Meanwhile, the aforementioned districts 7 and 11 barely landed above Category 1. Indeed, each of them seem reasonably compact with light noise. District 11 seems more irregular near its southeastern end, however.

\subsubsection{Category 1}

This category happily contains most of Florida's districts. Few visual issues are seen among these districts, aside from shared borders with districts of higher category discussed above. Indeed, district 25 has a jut from its southeast corner without much motivation, but each section of the district is rectangular without much thinness.

\subsection{Statewide analysis}

\begin{figure}
\centerline{\includegraphics[height=3.0in]{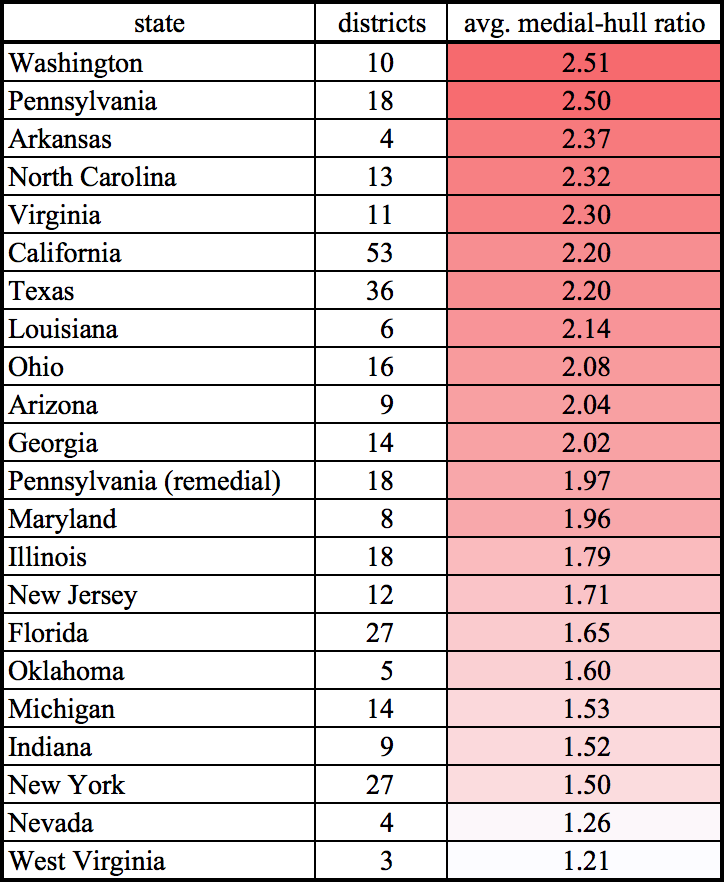}}
\caption{\label{fig:states} Statewide medial-hull ratio averages}
\end{figure}

Finally, we will average the medial-hull ratios of each\footnote{We completed calculations for all states with at least 10 Congressional districts, as well as for some with fewer. In total, we analyzed 335 districts from 22 states.} state's districts, then assign that value to the state. We may use these values to compare states and perhaps judge if a state's map overall tends to demonstrate evidence of gerrymandering.

At first glance over Figure \ref{fig:states}, the reader will notice a fairly even spread of values between 1.50 and 2.40. Indeed, we may consider applying the categories from Section \ref{sec:computation} to the states in aggregate rather than just to districts. With the category system, almost every state is either Category 1 or Category 2.

Notably, Pennsylvania is one of only two states assigned to Category 3 (and none were assigned to Category 4). The remedial map for Pennsylvania, however, earns it a reclassification into Category 1. We make special note of Washington, which tops the list: its districts average 2.01 except for two outliers, WA 1 and WA 8, which fall especial victim to the medial issue detailed in Section \ref{subsec:noise}. Though Washington is technically deemed Category 3, we will consider it as Category 2.

In this light, we may consider Category 3, and its suggestion of {\sl probable} gerrymandering, to call necessity to the redistricting of a state's Congressional map. We view categorical extremity as a measure of urgency. States high in Category 2, such as Arkansas, North Carolina, and Virginia, could certainly stand to be redistricted so that their average medial-hull ratios are lower, but they lack the need for immediate attention like that shown by Pennsylvania in early 2018. States like North Carolina currently have court cases actively pushing for this, and we ultimately seek in this paper to provide mathematical rigor as a tool to inform these legal decisions.

\section{Concluding Remarks}\label{sec:conclusion} No single mathematical measure can definitively identify gerrymandered districts. The medial-hull ratio developed here is another tool that can be used to analyze proposed maps for potential problems. As demonstrated above, this quantity often detects districts that have been ruled problematic by courts, and it is a reliable measure of how much a district wanders. In conjunction with other tools, the medial-hull ratio should prove useful to commissions tasked with drawing voting district maps.

\end{document}